\documentstyle[bezier]{article}
\begin{document}
\def\emline#1#2#3#4#5#6{%
       \put(#1,#2){\special{em:moveto}}%
       \put(#4,#5){\special{em:lineto}}}
\baselineskip=18pt

\begin{center}
\vspace{12mm}
{\Large \bf Partitions and their lattices}

\vspace{12mm}
Milan Kunz
 
{\bf Abstract}
\end{center}

Ferrers graphs and tables of partitions are treated as vectors.  
Matrix operations are used for simple proofs
of identities concerning partitions. Interpreting partitions
as vectors gives a possibility to generalize partitions on
negative numbers. Partitions are then tabulated
into lattices and some properties of these lattices are studied. 
There appears a new identity counting isoscele Ferrers graphs. 
The lattices form the base for tabulating combinatorial identities.

\section{Preliminary Notes} 
\label{Preliminary Notes}

Partitions of a natural number m into n parts were introduced into mathematics 
by Euler. The analytical formula for finding the number of partitions was 
derived by Ramanudjan and Hardy [1]. Ramanudjan was a mathematical 
genius from India. He was sure that it was possible to calculate the number of 
partitions exactly for any number m. He found the solution in cooperation with 
his tutor, the English mathematician Hardy. It is rather complicated formula 
derived by higher mathematical techniques. We will use only simple recursive 
methods for different relations between partitions.

Steve Weinberg in his lecture [2] about importance of mathematics for 
physics mentioned that partitions got importance for theoretical physics, even 
if Hardy did not want to study practical problems. But 
partitions were used in physics before Hardy by Boltzmann [3]. He
used this notion for splitting $m$ quanta of 
energy between $n$ particles in connection with his notion of entropy. 
He called partitions complexions, considering them to be orbits
in phase space. His idea was forgoten.

A partition splits a number $m$ into $n$ parts which sum is equal to the~number 
$m$, say $7:\ \ 3, 2, 1, 1$. A partition is an ordered set. Its objects, {\em 
parts}, are written in a row in decreasing order: 

$$ m_{j-1} \geq m_j \geq m_{j+1}\;.$$

If we close a string of parts into brackets, we get a $n$~dimensional vector 
row $p = (3, 2, 1, 1)$. From a partition vector, another vectors having 
equivalent structure of elements, for example. $r = (1, 2, 1, 3)$, are obtained 
by permuting, simple changing of ordering of vector elements. The 
partitions are thus indispensable for obtaining combinatorial 
identities, for  
ordering points of plane simplexes having constant sums of its constituting vectors. 

All unit permutations of a vector have the same length. Therefore different 
partitions form bases for other vectors composed from the same parts. Vectors 
belonging to the same partition of p into three parts are connected with 
other points of~the~three 
dimensional simplex by circles. In higher dimensions the circles become spheres 
and therefore we will call an ordered partition 
the {\em partition orbit} or simply orbit.

The number of vectors in partitions will be given as $n$, the size of the first 
vector as $m_1$. The bracket $(m,n)$ means all partitions of the number $m$ into 
at most $n$ parts. Because we write a partition as a n~dimensional vector we 
allow zero parts in a partition to fill empty places of the vector. It is a 
certain innovation against the tradition which will be very useful. But it is 
necessary to distinguish strictly both kinds of partitions, with zeroes and 
without them.

\section{Transposing and Transversing} 

Transposition of matrix vectors is the basic operation on matrices. 
It changes simply the row 
indices $i$ and column indices $j$ of all matrix elements

\begin{equation} 
{\bf M}^{\rm T} \rightarrow {\rm m}^{\rm T}_{ij} = {\rm m}_{ji}\; 
\end{equation}

\begin{figure} 
\caption{Transposing (A) and transversing (B) of matrices} 
\label{Transposing} 
\linethickness{0.6pt} 
\begin{picture}(140.00,90.00) 
\put(10.33,20.00){\framebox(39.67,40.00)[cc]{}} 
\put(79.67,20.00){\framebox(40.33,40.33)[cc]{}} 
\multiput(5.00,65.00)(0.12,-0.12){459}{\line(0,-1){0.12}} 
\multiput(69.67,10.00)(0.12,0.12){459}{\line(1,0){0.12}} 
\put(60.00,20.00){\vector(-1,2){0.2}} 
\bezier{184}(50.00,9.67)(70.33,-1.00)(60.00,20.00) 
\put(69.67,20.00){\vector(1,3){0.2}} 
\bezier{180}(79.67,9.67)(60.67,-1.00)(69.67,20.00) 
\put(30.00,70.00){\makebox(0,0)[cc]{A}} 
\put(100.00,70.00){\makebox(0,0)[cc]{B}} 
\end{picture} 
\end{figure}
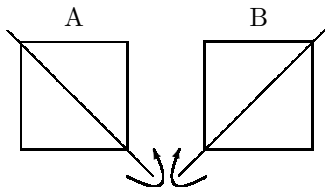

The second operation introduced here, {\em transvesing}, is not used in 
textbooks but we need it to prove simply, without calculations, some combinatorial 
identities concerning partitions. The transversing changes the ordering of both indices, that means 
rows and columns are counted backwards. If transposing rotates matrix elements 
around the main diagonal $m_{11} \longrightarrow m_{nn}$, transversing rotates 
them around the diagonal (its name will be transversal) $m_{1n} \longrightarrow 
m_{n1}$ (Fig. \ref{Transposing}). We look on the matrix's most distant corner as 
its starting point.

\section{Ferrers Graphs} 
\label{Ferrers}

\begin{figure} 
\caption{Ferrers graphs construction. New boxes are added to free places} 
\label{Ferrers Graphs} 
\linethickness{0.6pt} 
\begin{picture}(170.00,160.00) 
\put(70.00,133.00){\framebox(10.33,10.00)[cc]{}} 
\put(40.00,113.00){\framebox(10.00,10.00)[cc]{}} 
\put(40.00,103.00){\framebox(10.00,10.00)[cc]{}} 
\put(90.00,113.00){\framebox(10.00,10.00)[cc]{}} 
\put(100.00,113.00){\framebox(10.00,10.00)[cc]{}} 
\put(9.67,73.00){\framebox(10.33,10.00)[cc]{}} 
\put(9.67,83.00){\framebox(10.33,10.00)[cc]{}} 
\put(9.67,93.00){\framebox(10.33,10.00)[cc]{}} 
\put(60.33,83.00){\framebox(9.67,10.00)[cc]{}} 
\put(60.33,93.00){\framebox(9.67,10.00)[cc]{}} 
\put(70.00,93.00){\framebox(10.00,10.00)[cc]{}} 
\put(100.00,93.00){\framebox(10.00,10.00)[cc]{}} 
\put(110.00,93.00){\framebox(10.00,10.00)[cc]{}} 
\put(120.00,93.00){\framebox(10.00,10.00)[cc]{}} 
\put(9.67,33.33){\framebox(10.33,10.00)[cc]{}} 
\put(9.67,43.33){\framebox(10.33,10.00)[cc]{}} 
\put(9.67,53.33){\framebox(10.33,10.00)[cc]{}} 
\put(60.33,43.33){\framebox(9.67,10.00)[cc]{}} 
\put(60.33,53.33){\framebox(9.67,10.00)[cc]{}} 
\put(70.00,53.33){\framebox(10.00,10.00)[cc]{}} 
\put(111.00,30.33){\framebox(10.00,10.00)[cc]{}} 
\put(121.00,30.33){\framebox(10.00,10.00)[cc]{}} 
\put(131.00,30.33){\framebox(10.00,10.00)[cc]{}} 
\put(9.67,23.33){\framebox(10.33,10.00)[cc]{}} 
\put(70.00,43.33){\framebox(10.00,10.00)[cc]{}} 
\put(141.00,30.33){\framebox(10.00,10.00)[cc]{}} 
\put(27.67,33.00){\framebox(10.33,10.00)[cc]{}} 
\put(27.67,43.00){\framebox(10.33,10.00)[cc]{}} 
\put(27.67,53.00){\framebox(10.33,10.00)[cc]{}} 
\put(38.00,53.00){\framebox(10.00,10.00)[cc]{}} 
\put(93.33,53.00){\framebox(10.00,10.00)[cc]{}} 
\put(103.33,53.00){\framebox(10.00,10.00)[cc]{}} 
\put(113.33,53.00){\framebox(10.00,10.00)[cc]{}} 
\put(93.33,43.33){\framebox(10.00,9.67)[cc]{}} 
\put(50.00,123.00){\vector(-2,-1){0.2}} 
\multiput(69.67,133.00)(-0.23,-0.12){84}{\line(-1,0){0.23}} 
\put(20.00,103.00){\vector(-1,0){0.2}} 
\put(40.00,103.00){\line(-1,0){20.00}} 
\put(9.67,63.33){\vector(0,-1){0.2}} 
\put(9.67,73.00){\line(0,-1){9.67}} 
\put(38.00,63.00){\vector(2,-3){0.2}} 
\multiput(20.00,93.00)(0.12,-0.20){151}{\line(0,-1){0.20}} 
\put(100.00,123.00){\vector(2,-1){0.2}} 
\multiput(80.00,133.00)(0.24,-0.12){84}{\line(1,0){0.24}} 
\put(70.00,103.00){\vector(2,-1){0.2}} 
\multiput(50.00,113.00)(0.24,-0.12){84}{\line(1,0){0.24}} 
\put(120.00,103.00){\vector(1,-1){0.2}} 
\multiput(110.00,113.00)(0.12,-0.12){84}{\line(0,-1){0.12}} 
\put(80.00,103.00){\vector(-1,-1){0.2}} 
\multiput(90.33,113.00)(-0.12,-0.12){84}{\line(-1,0){0.12}} 
\put(38.00,43.00){\vector(-1,-2){0.2}} 
\multiput(60.33,83.00)(-0.12,-0.21){187}{\line(0,-1){0.21}} 
\put(70.00,53.33){\vector(0,-1){0.2}} 
\put(70.00,83.00){\line(0,-1){29.67}} 
\put(113.33,63.00){\vector(1,-1){0.2}} 
\multiput(80.00,93.00)(0.13,-0.12){251}{\line(1,0){0.13}} 
\put(141.00,40.33){\vector(1,-4){0.2}} 
\multiput(130.00,93.00)(0.12,-0.57){92}{\line(0,-1){0.57}} 
\put(93.33,52.67){\vector(-1,-4){0.2}} 
\multiput(100.00,93.00)(-0.12,-0.72){56}{\line(0,-1){0.72}} 
\end{picture} 
\end{figure}
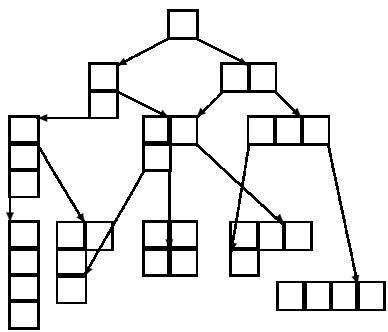

Ferrers graphs are used in the theory of partitions for many proofs based simply 
on their properties. Ferrers graphs are tables (see Fig. \ref{Ferrers Graphs}) 
containing $m$ objects, each object in its own box. The square boxes are 
arranged into columns in nonincreasing order $m_j \geq m_{j+1}$ with the sum

\begin{equation} 
\sum^n_{j=1} m_j = \sum^{\infty}_{k=0}n_km_k =m \;. 
\end{equation}

If partitions contain equal parts, it is possible to count them together using 
the index k and their number $n_k$.

It is obvious that a Ferrers graph, completed to an quadrangle with zero positions,
is a matrix ${\bf F}$ which has its unit 
elements arranged consecutively in the initial rows and columns. 

Introducing Ferrers graphs as matrices, we come necessarily to the notion of 
{\em restricted partitions}. The parts of a partitions can not be greater than 
the~number of rows of the matrix and the~number of parts greater than the~number 
of its columns.

\begin{figure} 
\caption{Truncation of partitions by restrictions of rows and columns} 
\label{Truncation} 
\linethickness{0.4pt} 
\begin{picture}(100.00,100.00) 
\put(20.00,19.67){\line(1,0){60.67}} 
\put(20.00,19.67){\line(0,1){60.33}} 
\multiput(20.00,70.00)(0.12,-0.12){420}{\line(0,-1){0.12}} 
\put(20.00,51.00){\line(1,0){18.67}} 
\put(47.33,42.67){\line(0,-1){23.00}} 
\put(20.00,53.00){\line(1,0){16.67}} 
\put(20.00,55.33){\line(1,0){14.33}} 
\put(20.00,57.33){\line(1,0){12.67}} 
\put(20.00,59.33){\line(1,0){10.67}} 
\put(20.00,61.00){\line(1,0){9.00}} 
\put(20.00,63.00){\line(1,0){7.00}} 
\put(20.00,65.00){\line(1,0){5.00}} 
\put(20.00,66.67){\line(1,0){3.33}} 
\put(20.00,68.00){\line(1,0){1.67}} 
\put(49.33,40.33){\line(0,-1){20.67}} 
\put(51.67,38.33){\line(0,-1){18.67}} 
\put(54.33,35.67){\line(0,-1){16.00}} 
\put(57.00,32.67){\line(0,-1){13.00}} 
\put(59.33,30.33){\line(0,-1){10.67}} 
\put(61.67,28.00){\line(0,-1){8.33}} 
\put(64.00,25.67){\line(0,-1){6.00}} 
\put(66.00,23.33){\line(0,-1){3.33}} 
\multiput(68.00,22.00)(0.11,-0.78){3}{\line(0,-1){0.78}} 
\put(10.00,69.67){\makebox(0,0)[cc]{m}} 
\put(9.67,51.00){\makebox(0,0)[cc]{M}} 
\put(47.33,10.00){\makebox(0,0)[cc]{N}} 
\put(70.33,10.00){\makebox(0,0)[cc]{n}} 
\end{picture} 
\end{figure}
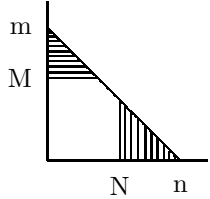

The interpretation of restrictions is geometrical. The part 
$m_{max}$ determines the side of a cube, $n$ is its dimension, see Fig. 
\ref{Truncation}.

A sophistication of the notation distinguishes the partitioned number M and the 
number of rows $m$ in the matrix  ${\bf F}$. The unrestricted number of partitions p(M) is equal to the 
number of restricted partitions when restricting conditions are loose then $m 
\geq M$ and $n \geq M$:

\begin{equation} 
p(M)_{unrestricted} = p(M,M,M)\;. 
\end{equation}

We write here first the number of rows $m$, then the number of parts $n$,
here considered as equal to $m$  and at
last the sum of unit elements (the number of filled boxes) M.

An important property of restricted partitions is determined by transposing 
Ferrers graphs ${\bf F} \rightarrow {\bf F}^{\rm T}$:

\begin{equation} 
p(m, n, M) = p(n, m, M)\;. 
\end{equation}

The partitions are conjugated. The number of partitions into exactly $n$~parts 
with the greatest part $m$ is the same as the number of partitions into $m$ 
parts having the greatest part $n$. This is simple transposing of ${\bf F}$.

A Ferrers graph can be subtracted from the matrix containing only unit elements 
(defined as ${\bf JJ}^{\rm T}$, {\bf J} being the unit column),
 and the resulting matrix transversed (Tr), for example

$$\begin{array}{ccccc} 
\left( 
\begin{array}{cc} 
1 & 1 \\ 
1 & 1 
\end{array} 
\right) & - & \left( 
\begin{array}{cc} 
1 & 0 \\ 
0 & 0 
\end{array} 
\right) & = & \left( 
\begin{array}{cc} 
0 & 1 \\ 
1 & 1 
\end{array} 
\right) 
\end{array}$$

$$\begin{array}{ccc} 
\left( 
\begin{array}{cc} 
0 & 1 \\ 
1 & 1 
\end{array} 
\right)^{\rm Tr} & = & \left( 
\begin{array}{cc} 
1 & 1 \\ 
1 & 0 
\end{array} 
\right) 
\end{array}$$

The relation between the number of restricted partitions of two different 
numbers is obtained according to the following equation

\begin{equation} 
p(m, n, M) = p(n, m, mn - M)\;. 
\end{equation}

This identity was derived by an operation very useful for acquiring elements of 
partition schemes (see later) and restricted partitions of all kinds. A 
restricted partition into exactly $n$ parts, having $m$ as the greatest part, 
has $(m + n - 1)$ units bounded by elements forming the first row and column of 
the corresponding Ferrers graph (Fig. \ref{Ferrers Graphs}). Only $(M - m - n + 
1)$ elements are free for partitions in the restricted frame $(m-1)$ and $(n-
1)$. Therefore

\begin{equation} 
p(m, n, M) = p(m-1, n-1, M-m-n+1)\;. 
\end{equation}

For example: p(4,3,8) = p(3,2,2) = 2. The corresponding partitions are 4,3,1 and 
4,2,2; or 2,0 and 1,1; respectively. This formula can be used for finding all 
restricted partitions. 

It is rather easy when the difference $(M-m-n+1)$ is 
smaller than the restricting values $m$ and $n$ or at least one from the 
restricting values. The row and column sums of partially 
restricted partitions 
having the~other constrain constant (shown as an asterics), 
where either $n$ or $m$ can be 1 till $M$ 
are:

\begin{equation} 
p(m, *, M) = \sum_{j=1}^M \ p(m, j, M) 
\end{equation}

\begin{equation} 
p(*, n, M) = \sum_{i=1}^M \ p(i, n, M)\;. 
\end{equation}

Before we examine restricted partitions in more detail, tables of unrestricted 
and partially restricted partitions will be introduced.

\section{Partition Matrices} 
\label{Partition Matrices}

Partially restricted partitions can be obtained from unrestricted partitions by 
subtracting a row of $n$ units or a column of $m$ units. This gives us 
the~recursive formula for the number of partitions as a sum of two partitions

\begin{equation} 
p(*, N, M) = p(*, N-1, M-1) + p(*, N, M-N-1)\;. 
\end{equation}

All partitions into exactly $N$ parts are divided into two sets. In one set are 
partitions having in the last column 1, their number is counted by the~term 
$p(*,N-1,M-1)$ which is the number of partitions of the number $(M-1)$ into 
exactly $(N-1)$ parts to which 1 was added on the n~th place and in other set 
are partitions which have in the last column 2 and more. They were obtained by 
adding the unit row ${\bf J}^{\rm T}$ with $n$ unit elements to the~partitions 
of $(M - N)$ into $N$ parts. Their number can be found in the same column
column $n$ places above.

A similar formula can be deduced for partitions of $M$ into at most $N$~parts. 
These partitions can have zero at least in the last column or they are 
partitioned into $n$ parts exactly:

\begin{equation} 
p(*, *=N, M) = p(*, *=N-1, M) + p(*, *=N, M-N)\;. 
\end{equation}

The term p(*, *=N-1, M) are partitions of $M$ into $(N - 1)$ parts transformed 
in partitions into $N$ parts by adding zero in the n-th column, the~term p(*, 
*=N, M-N) are partitions of $(M - 1)$ into $N$ parts to which the~unit row was 
added.

To formulate both recursive formulas more precisely, 
we had to define an apparently paradoxical partitionat first: 

$$p(0,0,0) = 1\;.$$

What it means? A partition of zero into zero number of parts. This partition 
represents the empty space of dimension zero. This partition is justified by its 
limit. We write $n = 0^0$ and find the limit:

\begin{equation} 
\lim{0^0} = \lim_{x\rightarrow\infty}\;(1/x)^0 = 1/x^0 = 1\;. 
\end{equation}

We get two following tables of partitions
  
\begin{table}
\caption{Partitions into exactly $n$ parts} 
\label{Partitions into Exactly $n$ Parts} 
\begin{tabular}{|r|rrrrrrr|r|} 
\hline 
n & 0 & 1 & 2 & 3 & 4 & 5 & 6 & $\Sigma$ \\ 
\hline 
m=0 & 1 & & & & & & & 1 \\ 
1 & & 1 & & & & & & 1 \\ 
2 & & 1 & 1 & & & & & 2 \\ 
3 & & 1 & 1 & 1 & & & & 3 \\ 
4 & & 1 & 2& 1 & 1& & & 5 \\ 
5 & & 1 & 2 & 2 & 1 & 1 & & 7 \\ 
6 & & 1 & 3 & 3 & 2 & 1 & 1 & 11 \\ 
\hline 
\end{tabular} 
\end{table}
  
\begin{table}
\caption{Partitions into at most n parts} 
\label{Partitions into at Most n Parts} 
\begin{tabular}{|r|rrrrrrr|} 
\hline 
n & 0 & 1 & 2 & 3 & 4 & 5 & 6 \\ 
\hline 
m=0 & 1 & 1 & 1 & 1 & 1 & 1 & 1 \\ 
1 & & 1 & 1 & 1 & 1 & 1 & 1 \\ 

2 & & 1 & 2 & 2 & 2 & 2 & 2 \\ 
3 & & 1 & 2 & 3 & 3 & 3 & 3 \\ 
4 & & 1 & 3 & 4 & 5 & 5 & 5 \\ 
5 & & 1 & 3 & 5 & 6 & 7 & 7 \\ 
6 & & 1 & 4 & 7 & 9 & 10& 11 \\ 
\hline 
\end{tabular} 
\end{table}

Table \ref{Partitions into at Most n Parts} is obtained from the 
Table \ref{Partitions into Exactly $n$ Parts} as partial sums of its rows, it means, 
by multiplying with the unit triangular matrix ${\bf T}^T$ from the right. The 
elements of the matrix ${\bf T}^{\rm T}$ are

\begin{equation} 
h_{ij} = 1\ {\rm if}\ j \geq i\;\ h_{ij} = 0 \ {\rm if}\ j > i\;. 
\end{equation}

On the other hand, the Table \ref{Partitions into at Most n Parts} 
is obtained from the Table \ref{Partitions into Exactly $n$ Parts} 
by multiplying 
with a matrix ${\bf T}^{\rm -T}$ from the right. The inverse elements are

\begin{equation} 
h_{ii}^{-1} = 1\;,\ h_{i,i+1}^{-1} = -1\;,\ h_{ij} = 0\;, {\rm otherwise}\;. 
\end{equation}

Notice, that the elements of the Table \ref{Partitions into at 
Most n Parts}  
right of the diagonal remain 
constant. They are equal to the row sums of the 
Table \ref{Partitions into Exactly $n$ Parts}. Increasing the number 
of zeroes does not change the number of partitions.

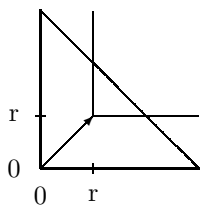
\begin{figure} 
\caption{Limiting of partition orbits. The lowest allowed part 
r shifts the plane simplex} 
\label{Limiting} 
\linethickness{0.4pt} 
\begin{picture}(100.00,100.00) 
\put(20.00,80.00){\line(0,-1){60.00}} 
\put(20.00,20.00){\line(1,0){60.00}} 
\put(40.00,79.33){\line(0,-1){39.33}} 
\put(40.00,40.00){\line(1,0){39.67}} 
\put(40.00,40.00){\vector(1,1){0.2}} 
\multiput(20.00,20.00)(0.12,0.12){167}{\line(0,1){0.12}} 
\multiput(20.00,80.00)(0.12,-0.12){501}{\line(0,-1){0.12}} 
\put(18.00,40.00){\line(1,0){4.00}} 
\put(40.00,22.00){\line(0,-1){4.00}} 
\put(10.00,40.00){\makebox(0,0)[cc]{r}} 
\put(10.00,20.00){\makebox(0,0)[cc]{0}} 
\put(40.00,10.00){\makebox(0,0)[cc]{r}} 
\put(20.00,10.00){\makebox(0,0)[cc]{0}} 
\end{picture} 
\end{figure}

When we multiply Table 
\ref{Partitions into Exactly $n$ Parts}  by the matrix ${\bf T}^{\rm -T}$ again, we obtain 
partitions having as the smallest allowed part the number 2. The effect of these 
operators can be visualized on the 2~dimensional complex, the operators shift 
the border of counted orbits (Fig. \ref{Limiting}). The operator ${\bf T}^{\rm 
T}$ differentiates n~dimensional complexes, shifting their border to positive 
numbers and cutting lover numbers. Zero forms the natural base border.

\section{Partitions with Negative Parts} 
\label{Partitions with Negative Parts}

Operations with tables of partitions lead to a thought, what would happen with 
partitions outside the positive cone of nonnegative numbers. Thus let us allow 
the existence of negative numbers in partitions, too\footnote{The negative parts 
can be compared in physics with antiparticles. Since an annihilation liberates 
energy, it does not annihilate it, the energy of the Universe is infinite. 
Speculations about existence of antiworlds, formed only by antiparticles 
balancing our world, can be formulated as doubt if the Universe is based in 
the~natural cone of the space.}.

If the number of equal parts $n_k$ is written as the vector row under the vector 
formed by the~number scale, the number of partitions is independent on shifts of 
the number scale, see Table \ref{Partitions as vectors}. Partitions, shown in the bottom 
part of the table,
 are always derived by shifting two 
vectors, one 1 position up, the other 1 position down. Each partition 
corresponds to a vector. If we write them as columns then their scalar product 
with the number scale, forming the vector row ${\bf m}^{\rm T}$, gives constant 
sum:

\begin{equation} 
{\bf m}^{\rm T}{\bf p} = \sum_{k \geq r}\; m_k n_k = m\;. 
\end{equation}

There is an inconsistency in notation, elements of the vector ${\bf p}$ are 
numbers of vectors having the same length and the letter $n$ with an index $k$ 
is used for them. For values of the number scale the letter $m$ is used with the 
common index $k$ which goes from the lowest allowed value of parts $r$ till the 
highest possible value. The index $k$ runs to infinity but all too high values 
$n_k$ are zeroes.
  
\begin{table}
\caption{Partitions as vectors} 
\label{Partitions as vectors} 
\begin{tabular}{|l|rrrrrr|lr|} 
\hline 
Parameter & r & & & & & & & \\ 
\hline 
Vector {\bf m} & -2 & -1 & 0 & 1 & 2 & 3 & ${\bf mp}$ = & -5 \\
 & -1& 0 & 1 & 2 & 3 & 4 & & 0 \\ 
& 0 & 1 & 2 & 3 & 4 & 5 & & 5 \\ 
& 1 & 2 & 3 & 4 & 5 & 6 & & 10 \\ 
& 2 & 3 & 4 & 5 & 6 & 7 && 15 \\ 
\hline 
Vector p & 4 & & & & & 1 & & \\ 
& 3 & 1 & & & 1 & & & \\ 
& 3 & & 1 & 1 & & & & \\ 
& 2 & 2 & & 1 & & & & \\ 
& 2 & 1 & 2 & & & & & \\ 
& 1 & 3 & 1 & & & & & \\ 
& 1 & 2 & 2 & & & & & \\ 
& & 5 & & & & & & \\ 
\hline 
\end{tabular} 
\end{table}

Using different partition vectors and different vectors ${\bf m}$ we get the 
following examples:

\begin{eqnarray*} 
(4 \times -2) + (1 \times 3) = -5\, \\ 
(3 \times -1) + (1 \times 0) + (1 \times 3) =\> 0\, \\ 
(3 \times 0) + (1 \times 2) + (1 \times 3) =\> 5\, \\ 
(2 \times 1) + (1 \times 2) + (2 \times 3) =\, 10\, \\ 
(1 \times 2) + (3 \times 3) + (1 \times 4) =\, 15. 
\end{eqnarray*}

The parameter $r$ shifts the table of partitions, its front rotates around 
the~zero point. If $r$ were $-\infty$, then $p(-\infty, 1) = 1$ but $p(-\infty, 
2)$ were undetermined, because a sum of a finite number with an infinite number 
is again infinite. The parameter $r$ will be written to a partition as its upper 
index to show that different bases of partitions are differentiating plane 
simplexes.

\section{Partitions with Inner Restrictions} 
\label{Partitions with Inner Restrictions}

Partitions were classified according to the minimal and maximal allowed values 
of parts, but there can be restrictions inside the number scale, it can be 
prescribed that some values are be forbidden. It is easy to see what this means: 
The plane simplex has holes, some orbits cannot be realized and its $(n -1)$ 
dimensional body is thinner than the normal one.

It is easy to find the number of partitions in which all parts are even. It is 
not possible to form an even partition from an uneven number, therefore:

\begin{equation} 
p_{\rm even}(2n) = p_{\rm unrestricted}(n)\;. 
\end{equation}

A more difficult task is finding the number of partitions in which all parts are 
odd. The rejected partitions contain mixed odd and even parts. The~relation 
between different partitions is etermined as

\begin{equation} 
p_{\rm unrestricted}(n) = p_{\rm odd}(n) + p_{\rm even}(n) + p_{\rm mixed}(n)\;. 
\end{equation}

The corresponding lists are given in Table 
\ref{Odd, even, and mixed partitions}

\begin{table}
\caption{Odd, even, and mixed partitions} 
\label{Odd, even, and mixed partitions} 
\begin{tabular}{|r|rrrrrrrrr|r|r|r|r|} 
\hline 
& \multicolumn{9}{|c|}{Number of odd partitions}& \multicolumn{4} {|c|}{ Sums}\\ 
\hline 
n & 1 & 2 & 3 & 4 & 5 & 6 & 7 & 8 & 9 & Odd & Even & Mixed & p(m) \\ 
\hline 
m=1 & 1 & & & & & & & & & 1 & 0 & 0 & 1 \\ 
2 & & 1 & & & & & & & & 1 & 1 & 0 & 2 \\ 
3 & 1 & & 1 & & & & & & & 2 & 0 & 1 & 3 \\ 
4 & & 1 & & 1 & & & & & & 2 & 2 & 1 & 5 \\ 
5 & 1 & & 1 & & 1 & & & & & 3 & 0 & 4 & 7 \\ 
6 & & 2 & & 1 & & 1 & & & & 4 & 3 & 4 & 11 \\ 
7 & 1 & & 2 & & 1 & & 1 & & & 5 & 0 & 10 & 15 \\ 
8 & & 2 & & 2 & & 1 & & 1 & & 6 & 5 & 11 & 22 \\ 
9 & 1 & & 3 & & 2 & & 1 & & 1 & 8 & 0 & 22 & 30\\ 
\hline 
\end{tabular} 
\end{table}

Notice how the scarce matrix of odd partitions is made from 
Table \ref{Partitions into Exactly $n$ Parts}. Its 
elements, except the first one in each column, are shifted down on cross 
diagonals. An odd number must be partitioned into an odd number of odd parts and 
an even number into even number of odd parts. Therefore the~matrix can be filled 
only in half. The recurrence is given by two possibilities how to increase the 
number $m$. Either we add odd 1 to odd partitions of~$(m - 1)$ with exactly $(j 
- 1)$ parts or we add $2j$ to odd numbers of partitions of $(m - 2j)$ with 
exactly $j$ parts. The relation is expressed as

\begin{equation} 
o(i,j) = p[(i+j)/2,j]\;. 
\end{equation}

Partitions with all parts unequal are important, because their transposed 
Ferrers graphs have the greatest part odd, when the number of parts is odd, and 
even, when the number of parts is even. For example

$$\begin{tabular}{rrrr} 
10 & & & \\ 
& 9,1 & & \\ 
& 8,2 & & \\ 
& 7,3 & 7,2,1 & \\ 
& 6,3 & 6,3,1 & \\ 
& & 5,4,1 & \\ 
& & 5,3,2 & \\ 
& & & 4,3,2,1 
\end{tabular}$$

The partitions with unequal parts can be tabulated as in 
Table \ref{Partitions with unequal parts}. Notice that 
the difference of the even and odd columns partitions is mostly zeroes and only 
sometimes $\pm 1$. The importance of this phenomenon will be explained later. 
The number of partitions with unequal parts coincide with the partitions which 
all parts are odd.
  
\begin{table}
\caption{Partitions with unequal parts} 
\label{Partitions with unequal parts} 
\begin{tabular}{|r|rrrr|c|c|} 
\hline 
n & 1 & 2 & 3 & 4 & $\Sigma$ & Difference ($n_{odd} - n_{even})$ \\ 
\hline 
m=1 & 1 & & & & 1 & 1 \\ 
2 & 1 & & & & 1 & 1 \\ 
3 & 1 & 1 & & & 2 & 0 \\ 
4 & 1 & 1 & & & 2 & 0 \\ 
5 & 1 & 2 & & & 3 & -1 \\ 
6 & 1 & 2 & 1 & & 4 & 0 \\ 
7 & 1 & 3 & 1 & & 5 & -1 \\ 
8 & 1 & 3 & 2 & & 6 & 0 \\ 
9 & 1 & 4 & 3 & & 8 & 0 \\ 
10 & 1 & 4 & 4 & 1 & 10 & 0 \\ 
11 & 1 & 5 & 5 & 1 & 12 & 0 \\ 
12 & 1 & 5 & 7 & 2 & 15 & 1 \\ 
\hline 
\end{tabular} 
\end{table}

The differences are due to Franklin blocks with growing minimal parts and 
growing number of parts (their transposed notation is used), which are minimal 
in that sense that their parts differ by one, the shape of corresponding Ferrers 
graphs is trapeze:

$$\begin{array}{ccccc} 
\begin{array}{c} 
(1)\\ \\ 
\left( 
\begin{array}{ccc} 
1 & 1 & 1 \\ 
1 & 1 & 0 
\end{array} 
\right)\\ \\ 
\left( 
\begin{array}{ccccc} 
1 & 1 & 1 & 1 & 1 \\ 
1 & 1 & 1 & 1 & 0 \\ 
1 & 1 & 1 & 0 & 0 
\end{array} 
\right) 
\end{array} & \quad & \begin{array}{c} 
(1 1) \\ \\ 
\left( 
\begin{array}{cccc} 
1 & 1 & 1 & 1 \\ 
1 & 1 & 1 & 0 
\end{array} 
\right)\\ \\ 
\left( 
\begin{array}{cccccc} 
1 & 1 & 1 & 1 & 1 & 1 \\ 
1 & 1 & 1 & 1 & 1 & 0 \\ 
1 & 1 & 1 & 1 & 0 & 0 
\end{array} 
\right) 
\end{array} & \quad & \begin{array}{c} 
1, 2\\ \\ 
5, 7 \\ \\ \\ 
12, 15 \\ \\ \\ 
\end{array} 
\end{array}$$

\section{Differences According to Unit parts} 
\label{Differences according to unit parts}

We have arranged restricted partitions according to the number of nonzero parts 
in Table 1. It is possible to classify partitions according the number 
of~vectors in the partition having any value. Using value 1, we get another kind 
of partition differences as in 
Table \ref{Partitions According to Unit Parts} .
  
\begin{table}
\caption{Partitions Differentiated According to Unit Parts} 
\label{Partitions According to Unit Parts} 
\begin{tabular}{|r|rrrrrrr|} 
\hline 
n & 0 & 1 & 2 & 3 & 4 & 5 & 6 \\ 
\hline 
m=0 & 1 & & & & & & \\ 
1 & 0 & 1 & & & & & \\ 
2 & 1 & 0 & 1 & & & & \\ 
3 & 1 & 1 & 0 & 1 & & & \\ 
4 & 2 & 1 & 1 & 0 & 1 & & \\ 
5 & 2 & 2 & 1 & 1 & 0 & 1 & \\ 
6 & 4 & 2 & 2 & 1 & 1 & 0 & 1\\ 
\hline 
\end{tabular} 
\end{table}

The elements of the table are:

\begin{equation} 
p_{i0} = p(i) - p(i - 1),\ p_{ij} = p_{i-1,j-1}\ ,{\rm otherwise}\;. 
\end{equation}

Table \ref{Partitions According to Unit Parts} is obtained from the 
following Table \ref{Partitions and their Euler inversion} of rows 
of unrestricted 
partitions by multiplying it with the matrix ${\bf T}^{-1}$. The zero column of 
the Table \ref{Partitions According to Unit Parts} is the difference of two consecutive unrestricted partitions 
according to $m$. To all partitions of $p(m-k)$ were added k ones. The 
partitions in the~zero column contain only numbers greater than 1. These 
partitions can not be formed from lower partitions by adding ones and they are 
thus a difference of the partition function according to the number $n_1$. Since 
Table \ref{Partitions According to Unit Parts} is composed, it is the product of two matrices, its inverse is 
composed, too.

\section{Euler Inverse of Partitions} 
\label{Euler Inverse of Partitions}

If we write successive partitions as column or row vectors as in 
Table 7, 
which elements are

\begin{equation} 
p_{ij} = p(i - j + 1)\;, 
\end{equation}

we find rather easily its inverse matrix which is given in the second part 
of~the same Table.
  
\begin{table}
\caption{Partitions and their Euler inversion} 
\label{Partitions and their Euler inversion} 
\begin{tabular}{|r|rrrrrrc|rrrrrr|} 
\hline 
& \multicolumn{6}{|c}{Partition table}& \qquad & \multicolumn{6} 
{|c|}{Euler inversion} \\ 
\hline 
j & 0 & 1 & 2 & 3 & 4 & 5 & & 0 & 1 & 2 & 3 & 4 & 5 \\ 
\hline 
i=0 & 1 & & & & & & & 1 & & & & & \\ 
1 & 1 & 1 & & & & & & -1 & 1 & & & & \\ 
2 & 2 & 1 & 1 & & & & & -1 & -1 & 1 & & & \\ 
3 & 3 & 2 & 1 & 1 & & & & 0 & -1 & -1 & 1 & & \\ 
4 & 5 & 3 & 2 & 1 & 1 & & & 0 & 0 & -1 & -1 & 1 & \\ 
5 & 7 & 5 & 3 & 2 & 1 & 1 & & 1 & 0 & 0 & -1 & -1 & 1 \\ 
\hline 
\end{tabular} 
\end{table}

The nonzero elements in the first column of the Euler inversion (and similarly 
in the next columns which are only shifted down one row) appear at indices, 
which can be expressed by the Euler identity concerning the coefficients of 
expansion of

\begin{equation} 
(1 - t)(1 - t^2)(1 - t^3)... = 1 + \sum_{i = 1}^\infty\; (-1)^i\; [t^{3i^2 - 
i)/2} + t^{3i^2 + i)/2}]\;. 
\end{equation}

For example the last row of the partition Table \ref{Partitions and their Euler 
inversion} is eliminated by multiplying it with the Euler inversion as: 

$(7 \times 1) + (5 \times -1) +( 3 \times -1) + (2 \times 0) + (1 \times 0) + 
(1 \times 1) = 0$

when $i = 1$, there is the pair of indexes at t 1, 2 with negative sign; for $i 
= 2$ the pair is 5, 7; for $i = 3$ the pair is $-12, -15$ and so on. These 
numbers are the distances from the base partition. The inverse matrix becomes 
scarcer as $p(m)$ increases, as it was already shown in Franklin partitions 
above. All inverse elements are $-1,0,1$. The nonzero elements of the Euler 
polynomial are obtained as sums of the product

\begin{equation} 
\prod_{i=1}^\infty \; (1 - t^i)\;. 
\end{equation}

This is verified by multiplying several terms of the infinite product. If we 
multiply the Euler polynomial with its inverse function

\begin{equation} 
\prod_i=1^\infty (1 - t^i)^{-1}\;, 
\end{equation}

we obtain 1. From this relation follows that partitions are generated by 
the~inverse Euler function which is the {\em generating function} of partitions. 
Terms $t^i$ must be considered as representing unequal parts.

The Euler function has all parts $t^i$ different. We have 
constructed such 
partitions in Table \ref{Partitions with unequal parts}. If the coefficient at $t^i$ is obtained as the product 
of~the~even number of $(1 - t^i)$ terms then the sign is positive, and if it is 
the result of~the~uneven number of~terms then the sign is negative. The 
coefficients are determined by the difference of the number of partitions with 
odd and even number of~unequal parts. This difference can be further explained 
according to Franklin using Ferrers graphs.

All parts in $p(n)$ having as at least one part equal to 1 are obtained from 
$p(n-1)$. The difference $p(n) - p(n-1)$ is due to some terms of $p(n-2)$. We 
must add to each partition of $p(n-2)$ 2, except all partitions of $p(n-2)$ 
containing 1. These must be either removed or used in transposed form using 
transposed Ferrers graphs, since big parts are needed. One from the pair of 
conjugate partitions is superfluous. These unused partitions must be subtracted. 
For example for $p(8)$:

$$\begin{array}{cccccc} 
6; & 1^6; & \qquad & Formed:& 8; & 62; \\ 
\underline{51};& 21^4;& & & 53; & \\ 
42; & 2^21^2;& & & 44;& 2^4; \\ 
\underline{33};& 2^3; & & & 3^22; & \\ 
41^2;& \underline{31^3}; & & & 42^2; & \\ 
\underline{321}; & & & & & 
\end{array}$$

Leftovers (underlined above):

\begin{center} 
p(1) + 5: 51; \qquad p(3) + 3: 33; 321; 31 
\end{center}

are obtained by subtracting the largest part from corresponding 
partition. Two 
must be added to the subtracted part. We get p(8-5) and p(8-7) as 
the~corrections.

\section{Other Inverse Functions of Partitions} 
\label{Other Inverse Functions of Partitions}

We already met other tables of partitions which have inverses 
because they are 
in lower triangular form. The inverse to the Table 
\ref{Partitions into Exactly $n$ Parts} is Table 
\ref{Inverse matrix to partitions into n parts}.
  
\begin{table}
\caption{Inverse matrix to partitions into n parts} 
\label{Inverse matrix to partitions into n parts} 
\begin{tabular}{|r|rrrrrr|} 
\hline 
n & 1 & 2 & 3 & 4 & 5 & 6 \\ 
\hline 
m=1 & 1 & & & & & \\ 
2 & -1 & 1 & & & & \\ 
3 & 0& -1 & 1 & & & \\ 
4 & 1 & -1 & -1 & 1 & & \\ 
5 & 0 & 1 & -1 & -1 & 1 & \\ 
6 & 0 & 1 & 0 & -1 & -1 & 1 \\ 
\hline 
\end{tabular} 
\end{table}

The inverse to Table \ref{Partitions According to Unit Parts} 
is Table \ref{Inverse matrix of unit differences}.
  
\begin{table}
\caption{Inverse matrix of unit differences}
\label{Inverse matrix of unit differences} 
\begin{tabular}{|r|rrrrrr|} 
\hline 
n & 1 & 2 & 3 & 4 & 5 & 6 \\ 
\hline 
m=1 & 1 & & & & & \\ 
2 & 0 & 1 & & & & \\ 
3 & -1 & 0 & 1 & & & \\ 
4 & -1 & -1 & 0 & 1 & & \\ 
5 & -1 & -1 & -1 & 0 & 1 & \\ 
6 & 0 & -1 & -1 & -1 & 0 & 1 \\ 
\hline 
\end{tabular} 
\end{table}

Whereas the columns of the Table \ref{Inverse matrix to partitions into n parts} are irregular and elements of each column 
must be found separately, columns of the Table \ref{Inverse matrix of unit differences} repeat as they are only 
shifted in each column one row down, similarly as the elements of~their parent 
matrix are. They can be easily found by multiplying the matrix of~the~Euler 
function (Table \ref{Partitions and their Euler inversion}) 
by the matrix ${\bf T}$ from the left.

\section{Partition Orbits in m Dimensional Cubes} 
\label{Partition Orbits in m Dimensional Cubes}

Restricted partitions have a geometric interpretation: They are orbits of 
n~dimensional plane complices truncated into cubes with the sides $(m -1)$ as on 
Fig. 3.

We can count orbits even in cubes. It is a tedious task if some special 
techniques are not applied, since their number depends on the size of the~cube. 
For example for the 3~dimensional space we get orbits as in Table 
\ref{Orbits in 3 dimensional cubes} .
  
\begin{table}
\caption{Orbits in 3 dimensional cubes}
\label{Orbits in 3 dimensional cubes} 
\begin{tabular}{|r|c|c|c|c|} 
\hline 
Edge size & 0 & 1 & 2 & 3 \\ 
\hline 
m=0 & 000 & 000 & 000 & 000 \\ 
1 & & 100 & 100 & 100 \\ 
2 & & 110 & 200; 110 & 210; 110 \\ 
3 & & 111 & 210; 111 & 300; 210; 111 \\ 
4 & & & 220; 211 & 310; 220; 211 \\ 
5 & & & 221 & 320; 311; 221 \\ 
6 & & & 222 & 330; 321; 222 \\ 
7 & & & & 331; 322 \\ 
8 & & & & 332 \\ 
9 & & & & 333 \\ 
\hline 
\end{tabular} 
\end{table}

The Equation 3 can be applied for cubes. It shows their important property, 
they are symmetrical along the main diagonal, going from the center of the 
coordinates, the simplex $n^0$ to the most distant vertex of the cube in~which 
all $n$ coordinates are $(m-1)$. The diagonal of~the~cube is represented on 
Table \ref{Orbits in 3 dimensional cubes} by $k$ indices. Moreover,
 a cube is convex, therefore

\begin{equation} 
M \leq mn/2 \ {\rm then} \ p(m,n,M) \geq p(m,n,M-1) 
\end{equation}

and if

\begin{equation} 
M \geq mn/2 \ {\rm then} \ p(m,n,M) \leq p(m,n,M-1)\;. 
\end{equation}

Here we see the importance of restricted partitions. From Table 10, we find 
the recurrence, which is given by the fact that in a greater cube the~lesser 
cube is always present as its base. New orbits which are on its enlarged sides 
are added to it. But it is enough to know orbits of one enlarged side, because 
the other sides are formed by these orbits. The enlarged side of a~n~dimensional 
cube is a $(n - 1) $~dimensional cube. The recurrence relation for~partitions in 
cubes is thus

\begin{equation} 
p(m,n,M) = p(m-1,n,M) + p(m,n-1,M)\;. 
\end{equation}

This recurrence will be explained later more thoroughly.

\section{Generating Functions of Partitions in Cubes} 
\label{Generating functions of partitions in cubes}

The generating function of partitions is simply the generating function 
of~the~infinite cube in the Hilbert space, which sides have different meshes:

\begin{equation} 
{\rm Parts\ 1:}\; (1 + t_1^1 + t_1^2 + \dots\; t_1^{\infty}) 
\end{equation}

\begin{equation} 
{\rm Parts\ 2:}\; (1 + t_2^1 + t_2^2 + \dots\; t_2^{\infty}) 
\end{equation}

and so on till

\begin{equation} 
{\rm Parts}\ \infty:\; (1 + \dots\; t^1_{\infty})\;. 
\end{equation}

When the multiplication's for all parts are made and terms on consecutive plane 
simplexes counted, we get:

\begin{equation} 
1 + t_1^1 + [t_2^1 + t_1^2] + [t_3^1 + \dots\;. 
\end{equation}

The generating function of restricted partitions is obtained by canceling 
unwanted (restricted) parts. Sometimes the generating function is formulated in 
an inverse form. The infinite power series are replaced by the differences $(1 - 
t_k^{-1})$. This is possible if we consider $t$ to be only a dummy variable. For 
example, the generating function of the partitions with unequal unrepeated parts 
is given by the product

\begin{equation} 
u(t) = \prod_{k=1}^\infty\ (1 - t_k)\;. 
\end{equation}

The mesh of the partition space is regular, it covers all numbers. The number of 
partitions is obtained by recursive techniques. But it is a very complicated 
function, if it is expressed in one closed formula, as the Ramanudjan-Hardy 
function is. The partitions form a carcass of the space. We will be interested, 
how the mesh of partitions is filled into the space which all axes have unit 
mesh and which contains also vector strings.

\section{Partition Schemes} 
\label{Partition Schemes}

Multidimensional plane simplexes are complicated objects and it is necessary to 
find tools how to analyze them. To draw them is impossible, as it was mentioned, 
because their parts are layered in our 3 dimensional world over themselves.

We already classified orbits in plane simplexes according to the number $k$ of 
nonzero parts. This number shows the dimensionality of subsimplices, their 
vertices, edges, and (k-1)~dimensional bodies. Lately we introduced the~number 
of unit vectors as a tool differentiating the simplex. Now we arrange partitions 
as two dimensional tables. These tables will be called {\em partition schemes}.

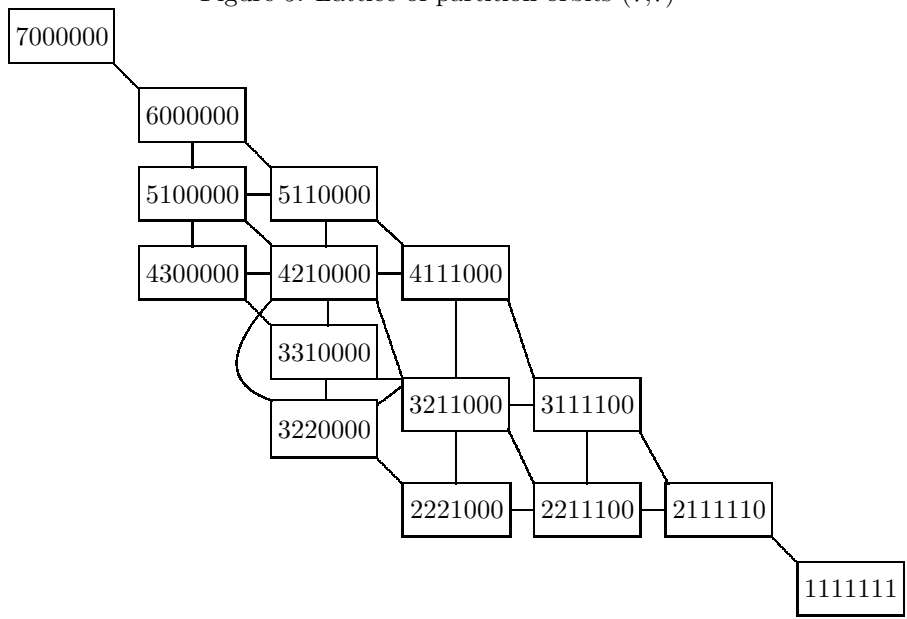
\begin{figure} 
\caption{Lattice of partition orbits (7,7)} 
\label{Lattice of Partition} 
\unitlength 0.7mm
\linethickness{0.4pt}
\begin{picture}(175.00,120.00)
\put(5.33,110.00){\framebox(19.67,10.00)[cc]{7000000}}
\put(30.00,95.00){\framebox(20.00,10.00)[cc]{6000000}}
\put(30.00,80.00){\framebox(20.00,10.00)[cc]{5100000}}
\put(30.00,65.00){\framebox(20.00,10.00)[cc]{4300000}}
\put(55.00,80.00){\framebox(20.00,10.00)[cc]{5110000}}
\put(55.00,65.00){\framebox(20.00,10.00)[cc]{4210000}}
\put(55.00,50.00){\framebox(20.00,10.00)[cc]{3310000}}
\put(55.00,35.00){\framebox(20.00,10.67)[cc]{3220000}}
\put(80.00,65.00){\framebox(20.00,10.00)[cc]{4111000}}
\put(80.00,40.00){\framebox(20.00,10.00)[cc]{3211000}}
\put(80.00,20.00){\framebox(20.33,10.00)[cc]{2221000}}
\put(105.00,40.00){\framebox(20.00,10.00)[cc]{3111100}}
\put(105.00,20.00){\framebox(20.00,10.00)[cc]{2211100}}
\put(130.00,20.00){\framebox(20.00,10.00)[cc]{2111110}}
\put(155.00,5.00){\framebox(20.00,10.00)[cc]{1111111}}
\multiput(25.00,110.00)(0.12,-0.12){42}{\line(0,-1){0.12}}
\put(40.00,95.00){\line(0,-1){5.00}}
\put(40.00,80.00){\line(0,-1){5.00}}
\multiput(50.00,95.00)(0.12,-0.12){42}{\line(1,0){0.12}}
\put(50.00,85.00){\line(1,0){5.00}}
\put(50.00,70.00){\line(1,0){5.00}}
\multiput(50.00,80.00)(0.13,-0.12){39}{\line(1,0){0.13}}
\multiput(50.00,65.00)(0.12,-0.12){42}{\line(1,0){0.12}}
\put(65.42,80.00){\line(0,-1){5.00}}
\put(65.83,65.00){\line(0,-1){5.00}}
\put(65.42,50.00){\line(0,-1){4.17}}
\bezier{136}(55.00,65.00)(41.67,50.42)(55.00,45.83)
\multiput(75.00,80.00)(0.13,-0.12){39}{\line(1,0){0.13}}
\put(75.00,70.00){\line(1,0){5.00}}
\put(90.00,65.00){\line(0,-1){15.00}}
\multiput(75.00,65.00)(0.12,-0.36){42}{\line(0,-1){0.36}}
\put(75.00,50.00){\line(1,0){5.00}}
\multiput(100.00,65.00)(0.12,-0.36){42}{\line(0,-1){0.36}}
\put(100.00,45.00){\line(1,0){5.00}}
\multiput(75.00,45.00)(0.16,0.12){32}{\line(1,0){0.16}}
\multiput(75.00,35.00)(0.12,-0.12){42}{\line(0,-1){0.12}}
\put(90.00,40.00){\line(0,-1){10.00}}
\multiput(125.00,40.00)(0.12,-0.22){46}{\line(0,-1){0.22}}
\put(115.00,40.00){\line(0,-1){10.00}}
\put(100.42,25.00){\line(1,0){4.58}}
\put(125.00,25.00){\line(1,0){5.00}}
\multiput(150.00,20.00)(0.12,-0.12){42}{\line(0,-1){0.12}}
\multiput(100.00,40.42)(0.12,-0.25){42}{\line(0,-1){0.25}}
\end{picture}
\end{figure}

Analyzing a 7~dimensional plane simplex with $m = 7$, we can start with its 3 
dimensional subsimplices. We see that they contain points 
corresponding to partitions: 7,0,0; 6,1,0; 5,2,0; 4,3,0; 5,1,1; 4,2,1; 3,3,1; 
3,2,2. The points corresponding to partitions are connected with other points 
of~the~simplex by circles. In higher dimensions the circles become spheres and 
this is the reason why we call a partition an {\em orbit}. The other points on 
each orbit have only different ordering of the same set of the coordinates.
  
\begin{table}
\caption{Partition scheme (7,7)} 
\label{Partition scheme (7,7)} 
\begin{tabular}{|r|rrrrrrr|r|} 
\hline 
n & 1 & 2 & 3 & 4 & 5 & 6& 7 & $\Sigma$ \\ 
\hline 
m = 7 & 1& & & & & & & 1 \\ 
6 & & 1 & & & & & & 1 \\ 
5 & & 1 & 1 & & & & & 2 \\ 
4 & & 1 & 1 & 1 & & & & 3 \\ 
3 & & & 2 & 1 & 1 & & & 4 \\ 
2 & & & & 1 & 1 & 1& & 3 \\ 
1 & & & & & & & 1 & 1 \\ 
\hline 
$\Sigma$ & 1 & 3 & 4 & 3 & 2 & 1 & 1 & 11 \\ 
\hline 
\end{tabular} 
\end{table}

Arranging partitions into tables (Table \ref{Partition scheme (7,7)}), the 
column classification is made according to the number of nonzero parts of 
partitions. Another classifying criterion is needed for rows. This will be the 
length of the longest vector $m_1$. From all partition vectors having the same 
dimensionality the longest vector is that one with the longest first vector. It 
dominates them. But there can exist longer orbits nearer to the surface of the 
simplex with a lesser number of nonzero parts. For example, vector (4,1,1) has 
equal length as (3,3,0) but vector (4,1,1,1,1) is shorter 
than (3,3,2,0,0). Such 
an arrangement is on Table \ref{Partition scheme (7,7)}. Orbits with three nonzero parts lie inside the 3 
dimensional simplex, with two nonzero parts lie on its edges. Orbits with four 
nonzero parts are inside tetrahedrons, it is on a surface in the fourth 
dimension. There exist these partitions: 4,1,1,1; 3,2,1,1; 2,2,2,1. Similarly 
columns corresponding to higher dimensions are filled.

The rows of partition schemes classify partitions according to the length of the 
first and longest vector ${\bf e}_1$. It can be shown easily that all vectors 
in~higher rows are longer than vectors in lover rows in corresponding columns. 
In the worst case it is given by the difference

\begin{equation} 
(x + 1)^2 + (x - 1)^2 > (2x)^2\;. 
\end{equation}

A three dimensional plane simplex to be can be considered as a truncated 
7~dimensional simplex, and after completing the columns of the Tab. 
\ref{Partition scheme (7,7)}) by the 
corresponding partitions, we get a crossection through the 7 
dimensional plane. 
The analysis is not perfect, an element is formed by two orbits, but 
nevertheless the scheme gives an insight how such high dimensional space looks 
like. We will study therefore properties of partitions schemes 
thoroughly.

The number of nonzero vectors in partitions will be given as $n$, the size 
of~the~first vector as $m$. Zeroes will not be written to spare work. The 
bracket $(m,n)$ means all partitions of the number $m$ into at most $n$ parts. 
Because we write a partition as a vector, we allow zero parts to complete the 
partition as before.

\section{Construction of Partition Schemes} 
\label{Construction of partition schemes}
  
\begin{table}
\caption{Partition scheme m = 13} 
\label{Partition Scheme m = 13} 
\begin{tabular}{|r|rrrrrr|rrrrrrr|} 
\hline 
n & 1 & 2 & 3 & 4 & 5 & 6 & 7 & 8 & 9 & 10 & 11& 12 & 13 \\ 
\hline 
m=13& 1 & & & & & & & & & & & & \\ 
12 & & 1& & & & & & & & & & & \\ 
11 & & 1& 1 & & & & & & & & & & \\ 
10 & & 1 & 1 & 1 & & & & & & & & & \\ 
9 & & 1 & 2 & 1 & 1 & & & & & & & & \\ 
8 & & 1 & 2 & 2 & 1 & 1 & & & & & & & \\ 
7 & & 1 & 3 & 3 & 2 & 1 & 1 & & & & & & \\ 
\hline 
6 & & & 3 & 4 & 3 & 2 & 1 & 1 & & & & & \\ 
5 & & & 2 & 4 & 5 & 3 & 2 & 1 & 1 & & & & \\ 
4 & & &   & 3 & 4 & 4 & 3 & 2 & 1 & 1 & & & \\ 
3 & & &     & & 2 & 3 & 3 & 2 & 2 & 1 & 1 & & \\ 
2 & & & & & &       & 1 & 1 & 1 & 1 & 1 & 1 & \\ 
1 & & & & & & & & & & & &  & 1 \\ 
\hline 
$\Sigma$ & 1 & 6 & 14 & 18 & 18 & 14 & 11 & 7 & 5 & 3& 2 & 1 & 1 \\ 
\hline 
\end{tabular} 
\end{table}

A partition scheme is divided into four blocks. Diagonal blocks repeat the~Table 
4.1 (the left upper block), the right lower one is written in the transposed 
form for $n > m/2$. Odd and even schemes behave somewhat differently, as can be 
seen on Tables \ref{Partition Scheme m = 13} and \ref{Partition Scheme m = 14}.
  
\begin{table}
\caption{Partition scheme m = 14} 
\label{Partition Scheme m = 14} 
\begin{tabular}{|r|rrrrrr|rrrrrrrr|} 
\hline 
n & 1 & 2 & 3 & 4 & 5 & 6 & 7 & 8 & 9 & 10 & 11& 12 & 13& 14 \\ 
\hline 
m=14& 1 & & & & & & & & & & & & & \\ 
13 & & 1 & & & & & & & & & & & & \\ 
12 & & 1 & 1 & & & & & & & & & & & \\ 
11 & & 1 & 1 & 1 & & & & & & & & & & \\ 
10 & & 1 & 2 & 1 & 1 & & & & & & & & & \\ 
9 & & 1 & 2 & 2 & 1 & 1 & & & & & & & & \\ 
8 & & 1 & 3 & 3 & 2 & 1& 1 & & & & & & & \\ 
7 & & 1 & 3 & 4 & 3 & 2 & 1 & 1 & & & & & & \\ 
\hline 
6 & & & 3 & * & * & * & 2 & 1 & 1 & & & & & \\ 
5 & & & 1 & * & * & * & 3 & 2 & 1 & 1 & & & & \\ 
4 & & & & 3 & * & * & 4 & 3 & 2 & 1 & 1& & & \\ 
3 & & & & & 2 & * & 3 & 3 & 2 & 2 & 1 & 1 & & \\ 
2 & & & & & & & 1 & 1 & 1 & 1 & 1 & 1 & 1 & \\ 
1 & & & & & & & & & & & & & & 1 \\ 
\hline 
$\Sigma$ & 1 & 7& 16 & 23 & 23 & 20 & 15 & 11 & 7 & 5 & 3 & 2 & 1 & 1\\ 
\hline 
\end{tabular} 
\end{table}

In the left lower block nonzero elements indicated by asterisks * can be placed 
only over the line which gives sufficient great product $mn$ to place all units 
into the corresponding Ferrers graphs and their sums must agree not only with 
row and column sums, but with diagonal sums, as we show below. This can be used 
for calculations of their numbers, together with rules for~restricted 
partitions.

The examples show three important properties of partition schemes:

\begin{itemize} 
\item Partition schemes are symmetrical according to their transversals, due to 
the conjugated partitions obtained by transposing Ferrers graphs. 
\item The upper left quarter (transposed lower right quarter) contain elements 
of the Table 4.1 of partitions into exactly $n$ parts shifted one column up. 
\item The schemes have form of the matrix in the lower diagonal form with unit 
diagonal. Therefore, they have inverses. It is easy to find them, for example 
for $n=7$ (Table \ref{Partition scheme (7,7) and its inversion}). 
\end{itemize}

The partitions in rows must be balanced by other ones with elements of~inverse 
columns. The third column includes or excludes 331 and 322 with 3211 and $31^4$; 
$2^31$ and $2^21^3$ with $2\times 21^5$, respectively.
  
\begin{table}
\caption{Partition scheme (7,7) and its inversion} 
\label{Partition scheme (7,7) and its inversion} 
\begin{tabular}{|r|rrrrrrr|crrrrrrr|} 
\hline 
n & 1 & 2 & 3 & 4 & 5 & 6 & 7 &\qquad &1 & 2 & 3 & 4 & 5 & 6 & 7 \\ 
\hline 
m = 7 & 1 & & & & & & & & 1 & & & & & & \\ 
6 & & 1 & & & & & & & & 1 & & & & & \\ 
5 & & 1 & 1 & & & & & & & 0 & 1 & & & & \\ 
4 & & 1 & 1 & 1 & & & & & & 0 & -1 & 1 & & & \\ 
3 & & & 2 & 1 & 1 & & & & & 2 & -1 & -1 & 1 & & \\ 
2 & & & & 1 & 1 & 1 & & & & -2 & 2 & 0 & -1& 1 & \\ 
1 & & & & & & & 1 & & & 0 & 0 & 0 & 0 & 0 & 1 \\ 
\hline 
\end{tabular} 
\end{table}

\section{Lattices of Orbits} 
\label{Lattices of orbits}

Partition orbit is a sphere which radius $r$ is determined by the Euclidean 
length of the corresponding vector: $r = (\sum p^2_j)$. Radiuses of some 
partition orbits coincide, for example $r(3,3,0)^2 = r(4,1,1)^2 = (18)$. It is 
thus impossible to determine distances between orbits using these radii 
(Euclidean distances) since the~distance between two different orbits cannot be 
zero.

We have shown in Sect. \ref{Partitions with Negative Parts} that one orbit 
can be obtained from another by 
shifting just two vectors, one up and other down on the number scale. We can 
imagine that both vectors collide and exchange their values as two particles of 
the ideal gas exchange their energy. If we limit the result of such an exchange 
to 1 unit, we can consider such two orbits to be the nearest neighbor orbits. 
The distance inside this pair is $\sqrt{2}$. We connect them in~the scheme by 
a~line. Some orbits are thus connected with many neighbor orbits, other have 
just one neighbor, compare with Fig. \ref{Lattice of Partition}. Orbits (3,3,0) and (4,1,1) are not 
nearest neighbors, because they must be transformed in~two steps: 

$$(3,3,0) \leftrightarrow ((3,2,1) \leftrightarrow (4,1,1)$$ 

or

$$(3,3,0) \leftrightarrow (4,2,0) \leftrightarrow (4,1,1)\;.$$

Partition schemes are generally not suitable for construction of orbit lattices, 
because at $m=n > 7$ there appear several orbits on some table places. It is 
necessary to construct at least 3~dimensional lattices to show all existing 
connections. For example:

$$\begin{array}{ccccc} 
(5,2,1) & \leftrightarrow & (4,3,1) & \leftrightarrow & (3,3,2) \\ 
& \searrow \nwarrow & \updownarrow & \swarrow \nearrow & \\ 
& & (4,2,2) & & 
\end{array}$$

\begin{figure} 
\caption{Lattice of file partitions. A file can be split into two new ones or 
two files can be combined into one} 
\label{Lattice of file partitions} 
\unitlength 0.50mm
\linethickness{0.6pt}
\begin{picture}(145.00,120.00)
\put(10.00,60.00){\framebox(10.00,10.00)[cc]{7}}
\put(30.00,60.00){\framebox(15.00,10.00)[cc]{52}}
\put(30.00,80.00){\framebox(15.00,10.00)[cc]{61}}
\put(30.00,40.00){\framebox(15.00,10.00)[cc]{43}}
\put(55.00,100.00){\framebox(15.00,10.00)[cc]{511}}
\put(55.00,70.00){\framebox(15.00,10.00)[cc]{421}}
\put(55.00,50.00){\framebox(15.00,10.00)[cc]{331}}
\put(55.00,20.00){\framebox(15.00,10.00)[cc]{322}}
\put(80.00,85.00){\framebox(20.00,10.00)[cc]{4111}}
\put(80.00,60.00){\framebox(20.00,10.00)[cc]{3211}}
\put(80.00,35.00){\framebox(20.00,10.00)[cc]{3211}}
\put(105.00,50.00){\framebox(15.00,10.00)[cc]{$2^21^3$}}
\put(125.00,60.00){\framebox(15.00,10.00)[cc]{$21^5$}}
\put(130.00,40.00){\framebox(15.00,10.00)[cc]{$1^7$}}
\multiput(20.00,70.00)(0.12,0.12){84}{\line(0,1){0.12}}
\put(20.00,70.00){\line(1,0){10.00}}
\multiput(20.00,60.00)(0.12,-0.12){84}{\line(0,-1){0.12}}
\multiput(45.00,90.00)(0.12,0.12){81}{\line(1,0){0.12}}
\multiput(70.00,99.67)(0.26,-0.12){39}{\line(1,0){0.26}}
\multiput(55.00,99.67)(-0.12,-0.35){84}{\line(0,-1){0.35}}
\multiput(45.00,80.00)(0.12,-0.24){84}{\line(0,-1){0.24}}
\multiput(45.00,60.00)(0.12,-0.36){84}{\line(0,-1){0.36}}
\put(45.00,50.00){\line(1,0){10.00}}
\multiput(45.00,50.00)(0.12,0.24){84}{\line(0,1){0.24}}
\put(45.00,70.00){\line(1,0){10.00}}
\multiput(45.00,40.00)(0.12,-0.12){84}{\line(0,-1){0.12}}
\multiput(70.00,99.67)(0.12,-0.35){84}{\line(0,-1){0.35}}
\multiput(70.00,80.00)(0.24,0.12){42}{\line(1,0){0.24}}
\multiput(70.00,70.00)(0.12,-0.30){84}{\line(0,-1){0.30}}
\put(70.00,60.00){\line(1,0){10.00}}
\multiput(70.00,50.00)(0.24,-0.12){42}{\line(1,0){0.24}}
\multiput(80.00,60.00)(-0.12,-0.36){84}{\line(0,-1){0.36}}
\multiput(70.00,30.00)(0.24,0.12){42}{\line(1,0){0.24}}
\put(100.00,85.00){\line(1,0){5.33}}
\multiput(100.00,85.00)(0.12,-0.52){48}{\line(0,-1){0.52}}
\multiput(100.00,70.00)(0.13,0.12){42}{\line(1,0){0.13}}
\put(100.00,60.00){\line(1,0){5.67}}
\multiput(100.00,45.00)(0.14,0.12){42}{\line(1,0){0.14}}
\put(120.00,60.00){\line(1,0){5.33}}
\multiput(125.33,60.00)(0.12,-0.26){39}{\line(0,-1){0.26}}
\put(45.00,80.00){\line(1,0){10.00}}
\put(70.00,70.00){\line(1,0){10.00}}
\put(105.67,74.67){\framebox(14.33,10.33)[cc]{$31^4$}}
\multiput(120.00,74.67)(0.14,-0.12){39}{\line(1,0){0.14}}
\end{picture}
\end{figure}
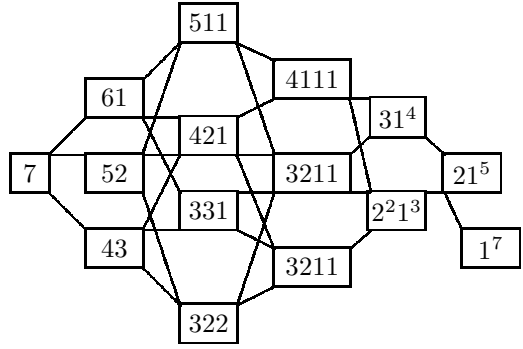

Sometimes stronger condition are given on processes going at exchanges, namely, 
that each collision must change the number of empty parts, as if they were 
information files which can be only joined into one file or one file separated 
into two or more files, or as if a part of a file transferred into an~empty 
file. Here also the nearest neighbor is limited on unifying of just 2 files or 
splitting a file into two (Fig.\ref{Lattice of file partitions}). In this case 
the path between two orbits must be longer, for example:

$$(3,3,0) \leftrightarrow (6,0,0) \leftrightarrow (4,2,0) \leftrightarrow 
(4,1,1)$$

or

$$(3,3,0) \leftrightarrow(3,2,1) \leftrightarrow(5,1,0) 
\leftrightarrow(4,1,1)\;.$$

In a lattice it is possible to count the number of nearest neighbors. If we 
investigate the number of one unit neighbors or connecting lines between columns 
of partition schemes, we obtain an interesting Table \ref{Right Hand One-unit 
Neighbors of Partition Orbits}.
  
\begin{table}
\caption{Right hand One-unit Neighbors of Partition Orbits} 
\label{Right Hand One-unit Neighbors of Partition Orbits}  
\begin{tabular}{|r|rrrrrr|r|} 
\hline 
n & 1 & 2 & 3 & 4 & 5 & 6 & $\Sigma$ \\ 
\hline 
m=2 & 1 & & & & & & 1 \\ 
3 & 1 & 1 & & & & & 2 \\ 
4 & 1 & 2 & 1 & & & & 4 \\ 
5 & 1 & 3 & 2 & 1 & & & 7 \\ 
6 & 1 & 4 & 4 & 2 & 1 & & 12 \\ 
7 & 1 & 5 & 6 & 4 & 2 & 1& 19 \\ 
\hline 
D(7-6) & 0 & 1 & 2 & 2 & 1& 1 & 7\\ 
\hline 
\end{tabular} 
\end{table}

The number of right hand neighbors is the sum of two terms. The right hand 
diagonal neighbors exist for all $p(m,n-1)$. We add 1 to all these partitions 
and decrease the largest part. Undetermined remain right hand neighbors in rows. 
Their number is equal to the number of partitions $p(m-2)$. To each partition 
$p(m-2, n-1)$ are added two units, one in the n~th column, the second in the (n-
1)~the column.

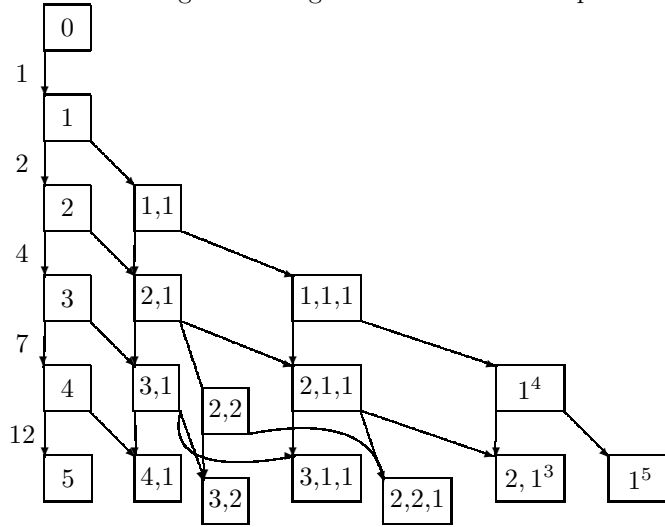
\begin{figure} 
\caption{Neighbor lattices between plane simplexes} 
\label{Neighbor lattices between plane simplexes} 
\unitlength 0.60mm
\begin{picture}(150.00,120.00)
\put(10.00,110.00){\framebox(10.00,10.00)[cc]{0}}
\put(10.00,90.00){\framebox(10.00,10.00)[cc]{1}}
\put(10.00,70.00){\framebox(10.00,10.00)[cc]{2}}
\put(30.00,70.00){\framebox(10.00,10.00)[cc]{1,1}}
\put(10.00,50.00){\framebox(10.00,10.00)[cc]{3}}
\put(30.00,50.00){\framebox(10.00,10.00)[cc]{2,1}}
\put(65.00,50.00){\framebox(15.00,10.00)[cc]{1,1,1}}
\put(10.00,30.00){\framebox(10.00,10.00)[cc]{4}}
\put(29.67,30.00){\framebox(10.00,10.00)[cc]{3,1}}
\put(45.00,25.00){\framebox(10.00,10.00)[cc]{2,2}}
\put(65.00,30.00){\framebox(15.00,10.00)[cc]{2,1,1}}
\put(110.00,30.00){\framebox(15.00,10.00)[cc]{$1^4$}}
\put(10.00,10.00){\framebox(10.33,10.00)[cc]{5}}
\put(30.00,10.00){\framebox(10.00,10.00)[cc]{4,1}}
\put(45.00,5.00){\framebox(10.00,10.00)[cc]{3,2}}
\put(65.00,10.00){\framebox(15.00,10.00)[cc]{3,1,1}}
\put(85.00,5.00){\framebox(15.00,10.00)[cc]{2,2,1}}
\put(110.00,10.00){\framebox(15.00,10.00)[cc]{$2,1^3$}}
\put(135.00,10.00){\framebox(14.00,10.00)[cc]{$1^5$}}
\put(10.00,100.00){\vector(0,-1){0.2}}
\put(10.00,110.00){\line(0,-1){10.00}}
\put(10.00,80.00){\vector(0,-1){0.2}}
\put(10.00,90.00){\line(0,-1){10.00}}
\put(10.00,60.00){\vector(0,-1){0.2}}
\put(10.00,70.00){\line(0,-1){10.00}}
\put(9.67,40.00){\vector(0,-1){0.2}}
\multiput(10.00,50.00)(-0.11,-3.33){3}{\line(0,-1){3.33}}
\put(10.00,20.00){\vector(0,-1){0.2}}
\put(10.00,30.00){\line(0,-1){10.00}}
\put(30.00,80.00){\vector(1,-1){0.2}}
\multiput(20.00,90.00)(0.12,-0.12){81}{\line(1,0){0.12}}
\put(65.33,60.00){\vector(3,-1){0.2}}
\multiput(40.00,70.00)(0.30,-0.12){84}{\line(1,0){0.30}}
\put(110.00,40.00){\vector(3,-1){0.2}}
\multiput(80.00,50.00)(0.36,-0.12){84}{\line(1,0){0.36}}
\put(135.00,20.00){\vector(1,-1){0.2}}
\multiput(125.00,30.00)(0.12,-0.12){84}{\line(1,0){0.12}}
\put(30.00,60.00){\vector(1,-1){0.2}}
\multiput(20.00,70.00)(0.12,-0.12){81}{\line(0,-1){0.12}}
\put(30.00,60.00){\vector(0,-1){0.2}}
\multiput(30.00,70.00)(-0.11,-3.33){3}{\line(0,-1){3.33}}
\put(30.00,40.00){\vector(1,-1){0.2}}
\multiput(20.00,50.00)(0.12,-0.12){81}{\line(0,-1){0.12}}
\put(30.00,20.00){\vector(1,-1){0.2}}
\multiput(20.00,30.00)(0.12,-0.12){84}{\line(0,-1){0.12}}
\put(30.00,40.00){\vector(0,-1){0.2}}
\put(30.00,50.00){\line(0,-1){10.00}}
\put(65.00,40.00){\vector(3,-1){0.2}}
\multiput(40.00,50.00)(0.30,-0.12){84}{\line(1,0){0.30}}
\multiput(40.00,50.00)(0.12,-0.36){42}{\line(0,-1){0.36}}
\put(65.00,40.00){\vector(0,-1){0.2}}
\put(65.00,50.00){\line(0,-1){10.00}}
\put(30.00,20.00){\vector(0,-1){0.2}}
\multiput(30.00,30.00)(0.11,-3.22){3}{\line(0,-1){3.22}}
\put(45.00,15.00){\vector(1,-3){0.2}}
\multiput(40.00,30.00)(0.12,-0.33){45}{\line(0,-1){0.33}}
\put(45.33,15.00){\vector(0,-1){0.2}}
\multiput(45.00,25.00)(0.11,-3.33){3}{\line(0,-1){3.33}}
\put(110.00,20.00){\vector(0,-1){0.2}}
\multiput(110.00,30.00)(-0.11,-3.33){3}{\line(0,-1){3.33}}
\put(109.67,20.00){\vector(3,-1){0.2}}
\multiput(80.33,30.00)(0.35,-0.12){84}{\line(1,0){0.35}}
\put(85.00,15.00){\vector(1,-3){0.2}}
\multiput(80.00,30.00)(0.12,-0.35){42}{\line(0,-1){0.35}}
\put(65.00,20.00){\vector(0,-1){0.2}}
\put(65.00,30.00){\line(0,-1){10.00}}
\put(65.33,20.00){\vector(1,0){0.2}}
\bezier{180}(40.00,30.00)(36.00,15.00)(65.33,20.00)
\put(85.00,15.00){\vector(1,-3){0.2}}
\bezier{168}(55.00,25.00)(82.00,30.00)(85.00,15.00)
\put(5.00,105.00){\makebox(0,0)[cc]{1}}
\put(5.00,85.00){\makebox(0,0)[cc]{2}}
\put(5.00,65.00){\makebox(0,0)[cc]{4}}
\put(5.00,45.00){\makebox(0,0)[cc]{7}}
\put(5.00,25.00){\makebox(0,0)[cc]{12}}
\end{picture} 
\end{figure}

The number of right hand neighbors $P(n)$ is the sum of the number 
of~unrestricted partitions

\begin{equation} 
P(n) = \sum_{k=0}^{n-2}\;p(k)\;. 
\end{equation}

To find all neighbors, we must add neighbors inside columns. The number of 
elements in columns is the number of partitions into exactly $n$ parts $p(m,n)$, 
the difference of each column must be decreased by 1 but there exist additional 
connections, see Fig. \ref{Neighbor lattices between plane simplexes}.

These connections must be counted separately. The resulting numbers are already 
known. The construction of partition schemes gives the result which we know as 
Table \ref{Partitions into Exactly $n$ Parts} read from the diagonal to the left.

The other interpretation of right hand one-unit neighbors of partitions is the 
plane complex as on Fig. \ref{Neighbor lattices between plane simplexes}. 
Vectors connect nearest neighbors in layers.

\section{Diagonal Differences in Lattices} 
\label{Diagonal Differences in Lattices}

In lattices, we can count orbits on side diagonals going consecutively parallel 
to the main diagonal. They count orbits having the form $[n - k]^k$. Their 
Ferrers graphs have a L form

$$\begin{array}{cccc} 
x & x & x & x \\ 
x & & & \\ 
x & & & \\ 
x & & & 
\end{array}$$

Side diagonal elements counts partitions which have in this layer smaller number 
of units, the other are inside this base.

The corresponding Table is \ref{Diagonal Sums of Partitions}.
  
\begin{table}
\caption{Diagonal Sums of Partitions} 
\label{Diagonal Sums of Partitions} 
\begin{tabular}{|r|rrrrrrrrr|r|} 
\hline 
k & 1 & 2 & 3 & 4 & 5 & 6 & 7 & 8 & 9 & $\Sigma$ \\ 
\hline 
n= 1 & 1 & & & & & & & & & 1 \\ 
2 & 2 & & & & & & & & & 2 \\ 
3 & 3 & & & & & & & & & 3 \\ 
4 & 4 & 1 & & & & & & & & 5 \\ 
5 & 5 & 2 & & & & & & & & 7 \\ 
6 & 6 & 3 & 2 & & & & & & & 11 \\ 
7 & 7 & 4 & 4 & & & & & & & 15 \\ 
8 & 8 & 5 & 6 & 3 & & & & & & 22 \\ 
9 & 9 & 6 & 8 & 6 & 1 & & & & & 30 \\ 
10 & 10 & 7 & 10 & 9 & 6 & & & & & 42 \\ 
11 & 11 & 8 & 12 & 12 & 11 & 2 & & & & 56 \\ 
12 & 12 & 9 & 14 & 15 & 16 & 9 & 2 & & & 77 \\ 
13 & 13 & 10 & 16 & 18 & 21 & 16 & 7 & & & 101 \\ 
14 & 14 & 11 & 18& 21& 26 & 23 & 18 & 4 & & 135 \\ 
15& 15 & 12& 20& 24& 31 & 30 & 29 & 12 & 3 & 176 \\ 
\hline 
\end{tabular} 
\end{table}

The initial $k$ column values have these analytical forms: 
\begin{itemize} 
\item 1n counts elements in n columns (rows) having the form $(n - k)1^k$, $k = 
0$ -- $(n -1)$; 
\item 1(n-3) counts elements in (n - 2) columns (rows) obtained from the basic 
partition 2,2 by adding units in the first row and column; 
\item 2(n-5) counts elements in (n - 2) columns (rows) obtained from the basic 
partitions 3,3 and 2,2,2 by adding units in the first row and column; 
\item 3(n-7) counts elements in (n - 2) columns (rows) obtained from the basic 
partitions 4,4; 3,3,2, and 2,2,2,2 by adding units in the first row and column; 
\item 5(n-9) + 1. On this level appears the partition 3,3,3 where elements start 
to occupy the third L layer; 
\item 7(n-11) + 2. 
\end{itemize}

The values in brackets are the numbers for partitions which lie inside the L 
frame having $(2k-1)$ units. At higher diagonal layers appear these 
possibilities to add new elements later. Partitions 4, 4, 4 
and 3, 3, 3, 3, for 
$n=12$, are counted in the seventh layer. For $n=13$, the layer counts seven 
partitions:

$$\begin{array}{ccccc} 
5,5,3; & & \\ 
5,4,4; & & \\ 
& 4,4,4,1; & \\ 
& 4,4,3,2; & \\ 
& 4,3,3,3; & \\ 
& & 3,3,3,3,1; \\ 
& & 3,3,3,2,1. 
\end{array}$$

There appears a very interesting property of partition lattices. 
The side diagonals being on side diagonals of the Table 
\ref{Diagonal Sums of Partitions} have equal length $n$, 
and the number of partitions $p(d)$ lying on them is equal to

\begin{equation} 
p(d)= 2^{(n-1)}
\end{equation} 

this is true for all complete diagonals in the table, also the seventh
diagonal sum is completed by the partition (4,4,4,4).
It can be conjectured, that it is a general property of lattices.
There are counted partitions which superposed Ferrers graphs can be
situated into  
isoscele triangular form ($M=N$) ending in the transversal
 which were not
counted before. The condition is that all Ferrers graphs are
superposed from the same starting place, otherwise Ferrers graphs 
of each partition can fill their isoscele triangular form.
 
The partitions can be ordered in the following way 
(see Table \ref{Binomial Ordering of Partitions}\.)

\begin{table}
\caption{Binomial Ordering of Partitions} 
\label{Binomial Ordering of Partitions} 
\begin{tabular}{|r|r|r|r|r|r|} 
\hline 
 1& 2 &3 &4 &5 & $\Sigma$ \\
\hline
(1)&  &   &   &  &  1 \\
\hline 
(1,1) & (2) & & & & 2 \\
\hline 
(1,1,1) & (2,1) & (3)  & & & \\
        & (2,2) & & & & 4 \\
\hline
(1,1,1,1) &(2,1,1)  & (3,1) & (4) & & \\ 
 & (2,2,1) &(3,2) & & & \\
 &(2,2,2)  &(3,3) & & & 8 \\
\hline
(1,1,1,1,1) & (2,1,1,1)  & (3,1,1) &(4,1) & (5) & \\ 
 &(2,2,1,1)  &(3,2,1) & (4,2) & & \\
 & (2,2,2,1) &(3,2,2) & (4,3)  & &  \\
 & (2,2,2,2) &(3,3,1) &(4,4) & & \\
    &       &(3,3,2) &             & &  \\
  &      &(3,3,3) &             & & 16  \\
\hline 
\end{tabular} 
\end{table}

Counting of partitions is changed into a combinatorial
problem of finding of all ordered combinations of
$(k-1)$ numbers with the greatest part equal to $k$.
The partitions are formed by a fixed part which is
equal to the number of the column and starts in the 
corresponding row. To this fixed part are added two
movable parts from the previous partitions, the whole upper
predecessor and the movable part of the left upper
predecessor. The resulting counts are the binomial numbers. 

\section{Generalized Lattices} 
\label{Generalized lattices}

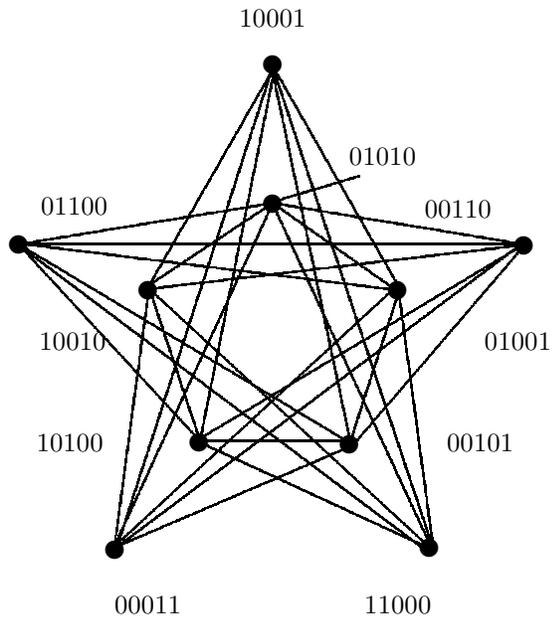
\begin{figure} 
\caption{Nearest neighbors in 00111 lattice} 
\label{Nearest Neighbors in 00111 Lattice} 
\unitlength 0.60mm
\linethickness{0.5pt} 
\begin{picture}(140.00,160.00) 
\put(70.00,129.67){\circle*{4.00}} 
\put(70.00,99.00){\circle*{4.00}} 
\put(13.67,90.00){\circle*{4.00}} 
\put(125.67,89.67){\circle*{4.00}} 
\put(42.33,79.67){\circle*{4.00}} 
\put(97.67,79.67){\circle*{4.00}} 
\put(53.67,46.00){\circle*{4.00}} 
\put(87.00,45.67){\circle*{4.00}} 
\put(104.67,22.67){\circle*{4.00}} 
\put(35.00,22.33){\circle*{4.00}} 
\put(70.00,140.00){\makebox(0,0)[cc]{10001}} 
\put(87.00,109.33){\makebox(0,0)[lc]{01010}} 
\put(18.67,98.33){\makebox(0,0)[lc]{01100}} 
\put(118.33,97.67){\makebox(0,0)[rc]{00110}} 
\put(25.67,68.33){\makebox(0,0)[cc]{10010}} 
\put(117.00,68.33){\makebox(0,0)[lc]{01001}} 
\put(32.33,46.00){\makebox(0,0)[rc]{10100}} 
\put(108.67,46.00){\makebox(0,0)[lc]{00101}} 
\put(35.00,10.00){\makebox(0,0)[lc]{00011}} 
\put(105.00,10.00){\makebox(0,0)[rc]{11000}} 
\multiput(70.33,129.33)(-0.12,-0.21){234}{\line(0,-1){0.21}} 
\multiput(70.33,129.00)(0.12,-0.21){231}{\line(0,-1){0.21}} 
\multiput(70.33,128.67)(-0.12,-0.60){139}{\line(0,-1){0.60}} 
\multiput(70.00,128.33)(0.12,-0.57){145}{\line(0,-1){0.57}} 
\multiput(70.00,129.33)(-0.12,-0.37){289}{\line(0,-1){0.37}} 
\multiput(70.00,129.33)(0.12,-0.36){292}{\line(0,-1){0.36}} 
\put(105.00,23.00){\line(0,-1){0.33}} 
\multiput(14.00,90.00)(0.74,0.12){76}{\line(1,0){0.74}} 
\multiput(13.33,90.33)(0.12,-0.13){337}{\line(0,-1){0.13}} 
\multiput(14.00,89.67)(0.20,-0.12){367}{\line(1,0){0.20}}
\multiput(14.00,90.00)(0.97,-0.12){87}{\line(1,0){0.97}} 
\put(13.67,90.00){\line(1,0){112.00}} 
\multiput(13.67,89.67)(0.16,-0.12){559}{\line(1,0){0.16}} 
\multiput(35.00,22.67)(0.12,0.90){64}{\line(0,1){0.90}} 
\multiput(35.00,23.00)(0.12,0.26){292}{\line(0,1){0.26}} 
\multiput(34.67,22.67)(0.13,0.12){473}{\line(1,0){0.13}} 
\multiput(35.67,23.00)(0.28,0.12){189}{\line(1,0){0.28}} 
\multiput(35.33,22.67)(0.16,0.12){562}{\line(1,0){0.16}} 
\multiput(53.67,46.00)(0.26,-0.12){198}{\line(1,0){0.26}} 
\multiput(43.00,79.67)(0.13,-0.12){473}{\line(1,0){0.13}} 
\multiput(70.00,99.00)(0.12,-0.26){292}{\line(0,-1){0.26}} 
\multiput(98.00,79.33)(0.12,-0.89){64}{\line(0,-1){0.89}} 
\multiput(126.00,90.00)(-0.12,-0.14){323}{\line(0,-1){0.14}} 
\multiput(126.00,90.00)(-0.20,-0.12){364}{\line(-1,0){0.20}} 
\multiput(125.33,89.67)(-1.02,-0.12){81}{\line(-1,0){1.02}} 
\multiput(125.33,90.00)(-0.73,0.12){76}{\line(-1,0){0.73}} 
\multiput(70.00,99.00)(-0.18,-0.12){156}{\line(-1,0){0.18}} 
\multiput(42.67,80.33)(0.12,-0.36){95}{\line(0,-1){0.36}} 
\put(54.00,46.33){\line(1,0){33.33}} 
\multiput(87.33,46.33)(0.12,0.37){89}{\line(0,1){0.37}} 
\multiput(98.00,79.33)(-0.17,0.12){164}{\line(-1,0){0.17}} 
\multiput(70.00,99.33)(0.40,0.12){48}{\line(1,0){0.40}} 
\end{picture} 
\end{figure}

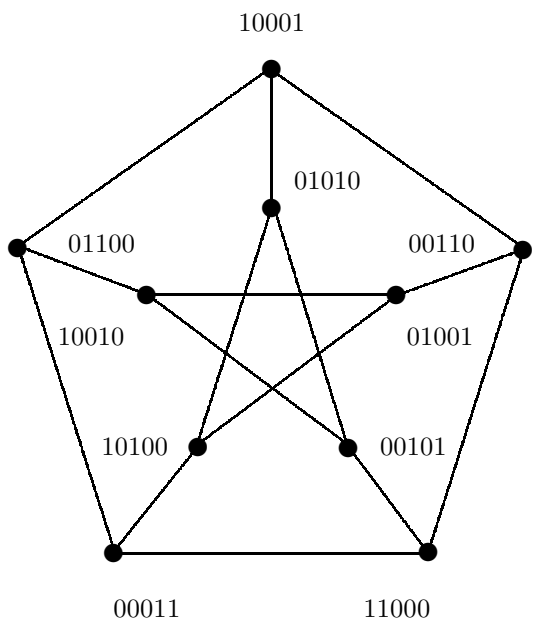
\begin{figure} 
\caption{Petersen graph. Adjacent vertices are in distances 4} 
\label{Petersen graph} 
\unitlength 0.60mm
\linethickness{0.5pt} 
\begin{picture}(140.00,160.00) 
\put(70.00,130.00){\line(0,-1){30.33}}
\multiput(70.00,130.00)(-0.17,-0.12){334}{\line(-1,0){0.17}} 
\multiput(14.00,90.00)(0.12,-0.38){176}{\line(0,-1){0.38}} 
\put(35.00,22.33){\line(1,0){69.67}} 
\multiput(104.67,22.33)(0.12,0.38){178}{\line(0,1){0.38}} 
\multiput(126.00,90.00)(-0.17,0.12){334}{\line(-1,0){0.17}} 
\multiput(104.67,22.33)(-0.12,0.16){148}{\line(0,1){0.16}} 
\multiput(87.00,45.33)(-0.12,0.39){139}{\line(0,1){0.39}} 
\multiput(70.33,100.00)(-0.12,-0.38){142}{\line(0,-1){0.38}} 
\multiput(53.33,45.67)(-0.12,-0.15){153}{\line(0,-1){0.15}} 
\multiput(14.00,90.33)(0.32,-0.12){87}{\line(1,0){0.32}} 
\multiput(42.00,80.00)(0.16,-0.12){289}{\line(1,0){0.16}} 
\put(42.00,79.67){\line(1,0){56.00}} 
\multiput(98.00,79.67)(0.33,0.12){87}{\line(1,0){0.33}} 
\multiput(98.00,79.67)(-0.16,-0.12){284}{\line(-1,0){0.16}} 
\put(70.00,129.67){\circle*{4.00}} 
\put(70.00,99.00){\circle*{4.00}} 
\put(13.67,90.00){\circle*{4.00}} 
\put(125.67,89.67){\circle*{4.00}} 
\put(42.33,79.67){\circle*{4.00}} 
\put(97.67,79.67){\circle*{4.00}} 
\put(53.67,46.00){\circle*{4.00}} 
\put(87.00,45.67){\circle*{4.00}} 
\put(104.67,22.67){\circle*{4.00}} 
\put(35.00,22.33){\circle*{4.00}} 
\put(70.00,140.00){\makebox(0,0)[cc]{10001}} 
\put(75.00,105.00){\makebox(0,0)[lc]{01010}} 
\put(25.00,91.00){\makebox(0,0)[lc]{01100}} 
\put(115.00,91.00){\makebox(0,0)[rc]{00110}} 
\put(30.00,70.33){\makebox(0,0)[cc]{10010}} 
\put(100.00,70.33){\makebox(0,0)[lc]{01001}} 
\put(32.33,46.00){\makebox(0,0)[lc]{10100}} 
\put(108.67,46.00){\makebox(0,0)[rc]{00101}} 
\put(35.00,10.00){\makebox(0,0)[lc]{00011}} 
\put(105.00,10.00){\makebox(0,0)[rc]{11000}} 
\end{picture} 
\end{figure}

The notion of lattices can be used also for possible transformations of points 
having specific properties among themselves, for example between all 10 
permutations of a 5~tuple composed from 3 symbols of one kind and 2 symbols of 
another kind. When the neighbors differ only by one exchange of the position of 
only anyone pair of two kinds symbols we obtain lattice as on 
Fig.\ref{Nearest Neighbors in 00111 Lattice}. Each from three unit symbols has two possibilities 
to change 0 into 1. Rearrange these ten points as a simple triangle. The 
simultaneous exchange of two pairs (or two consecutive changes of one pair give 
a pattern as on Fig.\ref{Petersen graph}, known as the Pettersen graph.

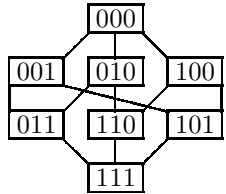
\begin{figure} 
\caption{Lattice of the three dimensional unit cube} 
\label{Lattice of the three dimensional unit cube} 
\linethickness{0.6pt} 
\begin{picture}(110.00,100.00) 
\put(40.33,70.00){\framebox(20.00,10.00)[cc]{000}} 
\put(40.33,50.00){\framebox(20.00,10.00)[cc]{010}} 
\put(10.33,50.00){\framebox(20.00,10.00)[cc]{001}} 
\put(70.33,50.00){\framebox(20.00,10.00)[cc]{100}} 
\put(10.33,30.00){\framebox(20.00,10.00)[cc]{011}} 
\put(40.33,30.00){\framebox(20.00,10.00)[cc]{110}} 
\put(70.33,30.00){\framebox(20.00,10.00)[cc]{101}} 
\put(40.33,9.67){\framebox(20.00,10.33)[cc]{111}} 
\multiput(30.33,60.00)(0.12,0.12){84}{\line(0,1){0.12}} 
\put(50.00,70.00){\line(0,-1){10.00}} 
\multiput(60.33,70.00)(0.12,-0.12){84}{\line(1,0){0.12}} 
\put(10.33,50.00){\line(0,-1){10.00}} 
\multiput(30.33,30.00)(0.12,-0.12){84}{\line(0,-1){0.12}} 
\put(50.00,50.00){\line(0,-1){10.00}} 
\put(50.00,30.00){\line(0,-1){10.00}} 
\multiput(70.33,30.00)(-0.12,-0.12){84}{\line(-1,0){0.12}} 
\put(90.33,50.00){\line(0,-1){10.00}} 
\multiput(40.33,50.00)(-0.12,-0.12){84}{\line(0,-1){0.12}} 
\multiput(30.33,50.00)(0.48,-0.12){84}{\line(1,0){0.48}} 
\multiput(70.33,50.00)(-0.12,-0.12){84}{\line(-1,0){0.12}} 
\end{picture} 
\end{figure}

\begin{figure} 
\caption{Four dimensional cube projection. One 3 dimensional cube is twisted 
$45^0$} 
\label{Four dimensional cube} 
\linethickness{0.6pt} 
\begin{picture}(150.00,150.00) 
\put(10.33,10.00){\framebox(120.00,120.00)[cc]{}} 
\put(50.00,50.00){\framebox(40.00,40.00)[cc]{}} 
\multiput(70.00,119.67)(-0.12,-0.12){414}{\line(-1,0){0.12}} 
\multiput(20.00,70.00)(0.12,-0.12){417}{\line(0,-1){0.12}} 
\multiput(70.00,20.00)(0.12,0.12){417}{\line(1,0){0.12}} 
\multiput(120.33,70.00)(-0.12,0.12){417}{\line(-1,0){0.12}} 
\multiput(70.00,84.67)(-0.12,-0.12){123}{\line(-1,0){0.12}} 
\multiput(55.00,70.00)(0.12,-0.12){126}{\line(0,-1){0.12}} 
\multiput(70.00,55.00)(0.12,0.12){126}{\line(0,1){0.12}} 
\multiput(85.00,70.00)(-0.12,0.12){123}{\line(0,1){0.12}} 
\put(10.33,130.00){\circle*{4.00}} 
\put(130.33,130.00){\circle*{4.00}} 
\multiput(10.67,129.67)(0.71,-0.12){84}{\line(1,0){0.71}} 
\multiput(10.33,130.00)(0.12,-0.12){331}{\line(0,-1){0.12}} 
\multiput(130.33,130.00)(-0.12,-0.12){334}{\line(-1,0){0.12}} 
\multiput(130.33,130.00)(-0.12,-0.69){87}{\line(0,-1){0.69}} 
\multiput(90.00,50.00)(0.12,-0.12){334}{\line(1,0){0.12}} 
\multiput(130.33,10.00)(-0.72,0.12){84}{\line(-1,0){0.72}} 
\multiput(10.67,10.00)(0.12,0.12){328}{\line(0,1){0.12}} 
\put(10.33,10.00){\line(0,1){0.00}} 
\multiput(10.33,10.00)(0.12,0.74){81}{\line(0,1){0.74}} 
\put(70.00,120.00){\line(0,-1){35.00}} 
\put(85.00,70.00){\line(1,0){35.33}} 
\put(70.00,55.00){\line(0,-1){35.00}} 
\put(20.33,69.67){\line(1,0){34.67}} 
\put(55.00,69.67){\line(-1,0){34.67}} 
\put(20.33,70.00){\line(1,0){34.67}} 
\multiput(90.00,90.00)(-0.12,-0.47){42}{\line(0,-1){0.47}} 
\multiput(90.00,50.00)(-0.48,0.12){42}{\line(-1,0){0.48}} 
\multiput(50.00,50.00)(0.12,0.48){42}{\line(0,1){0.48}} 
\multiput(50.00,90.00)(0.48,-0.12){42}{\line(1,0){0.48}} 
\put(70.00,120.00){\circle*{4.00}} 
\put(50.00,90.00){\circle*{4.00}} 
\put(90.00,90.00){\circle*{4.00}} 
\put(70.00,84.67){\circle*{4.00}} 
\put(20.00,70.00){\circle*{4.00}} 
\put(54.67,70.00){\circle*{4.00}} 
\put(84.67,70.00){\circle*{4.00}} 
\put(120.33,70.00){\circle*{4.00}} 
\put(50.00,50.00){\circle*{4.00}} 
\put(70.00,55.00){\circle*{4.00}} 
\put(90.00,50.00){\circle*{4.00}} 
\put(70.00,20.00){\circle*{4.00}} 
\put(10.00,10.00){\circle*{4.00}} 
\put(130.00,10.00){\circle*{4.00}} 
\end{picture} 
\end{figure}
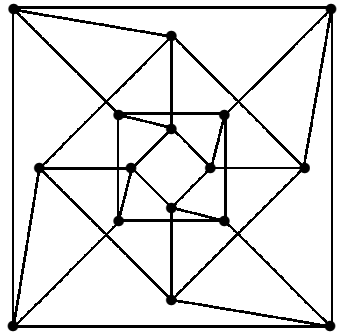

Lattices are formed by vertices of n~dimensional cubes. The nearest vertices 
differ only by one coordinate. The lattices of the 3~dimensional cube is on Fig. 
\ref{Lattice of the three dimensional unit cube}. Compare lines of the graphs 
with a real 3~dimensional cube and try to imagine the 4~dimensional cube (Fig. 
\ref{Four dimensional cube}).

A classical example of relation lattices is Aristotle's attribution of four 
properties: {\bf warm}, {\bf cold}, {\bf dry}, and {\bf humid} to four elements: 
fire, air, water and earth, respectively. It can be arranged in a form

$$\begin{array}{ccccc} 
air & & {\bf humid} & & water \\ \\ 
& {\bf warm} & 0 & {\bf cold} & \\ \\ 
fire & & {\bf dry} & & earth\;. 
\end{array}$$

The elements have always only two properties. The properties adjacent vertically 
and horizontally exclude themselves. Something can not be simultaneously warm 
and cold, or humid and dry\footnote { More precisely, it is necessary to draw a 
borderline (point zero) between these properties. Depending on its saturation, 
water vapor can be dry as well as wet.}.

\vspace{8mm}

{\bf{REFERENCES}}\\ 

[1] G. E. Andrews, The Theory of Partitions, Addison-Wesley Publ. 
Comp., Reading, MA, 1976.

[2] S. Weinberg, Mathematics, the Unifying Thread in Science, Notices 
AMS, 1986, 716.

[3] L. Boltzmann, \"Uber die Beziehung zwischen dem zweiten Hauptsatze 
der mechanishen W\"armetheorie und die Wahrscheinlichkeitsrechnung, {\it Wiener 
Berichte} {\bf 1877}, 76, 373.
\end{document}